\documentclass[%
  reprint,
  superscriptaddress,
  amsmath,amssymb,
  aps,
  floatfix,
  showkeys
]{revtex4-2}

\usepackage[
  activate={true,compatibility},
  final,
  tracking=true,
  kerning=true,
  nopatch=footnote
]{microtype}
\usepackage[T1]{fontenc}
\usepackage[utf8]{inputenc}

\usepackage[margin=1in]{geometry}

\usepackage{mathtools}
\usepackage{amsfonts}
\usepackage{amsmath}
\usepackage{nicefrac}
\usepackage{bbm}         
\usepackage{physics}
\usepackage{stmaryrd}

\usepackage{dsfont}
\usepackage{graphicx}
\usepackage{dcolumn}
\usepackage{booktabs}
\usepackage{tabularx}
\usepackage{multirow}
\usepackage{subcaption}
\usepackage{xcolor}
\usepackage[normalem]{ulem}
\usepackage{placeins}     %

\usepackage[noend]{algpseudocode}
\usepackage{algorithm}

\usepackage{amsthm}

\usepackage{hyperref}

\usepackage{tikz}
\usetikzlibrary{positioning}

\usepackage{chngcntr}
\counterwithout{equation}{section}  
\counterwithout{figure}{section}    

\makeatletter
\renewcommand\@makecaption[2]{%
  \par
  \vskip\abovecaptionskip
  \begingroup
    \small\rmfamily
    \begingroup
      \samepage
      \flushing
      \let\footnote\@footnotemark@gobble
      \@make@capt@title{#1}{#2}\par
    \endgroup
  \endgroup
  \vskip\belowcaptionskip
}
\makeatother

\newcolumntype{M}{>{\centering}X}
\newcolumntype{Y}{>{\hsize=.17\textwidth\arraybackslash}X}
\newcolumntype{C}{>{\hsize=.115\textwidth\centering\arraybackslash}X}
\newcolumntype{R}{>{\hsize=.115\textwidth\raggedleft\arraybackslash}X}

\begin{document}

\title{Efficient upsampling for tensor-network and quantum-state encoded functions}

\author{Siddhartha Guzman}
\affiliation{Quantum Research Center, Technology Innovation Institute, Abu Dhabi, UAE}

\author{Egor Tiunov}
\affiliation{Quantum Research Center, Technology Innovation Institute, Abu Dhabi, UAE}
\author{Leandro Aolita}
\affiliation{Quantum Research Center, Technology Innovation Institute, Abu Dhabi, UAE}
\date{\today}

\begin{abstract}
Both tensor trains (TTs) and quantum states provide compressed representations of grid-structured data with potentially exponential compression power. We present a unified framework for upsampling data encoded in vector amplitudes, with efficient realizations in both classical TT and quantum settings. Starting from an \(n\)-core TT or an \(n\)-qubit state on a coarse grid with \(2^n\) points, the construction produces an \((n+m)\)-core TT or \((n+m)\)-qubit state on a finer grid with \(2^{n+m}\) points. In the TT setting, it supports interpolation, quasi-interpolation, augmentation, and synthesis through efficient low-rank contractions, with the added \(m\) cores retaining constant rank. For function-value encodings, the resulting interpolation satisfies an \(\ell^2\)-error bound independent of the number of added grid points, achieves exponential compression at fixed accuracy, and has a logarithmic complexity in the number of grid points. In the quantum setting, the refined state is prepared by a \(\mathrm{poly}(n,m)\)-size circuit using \(\log(p+1)\) ancillas, where \(p\) controls the smoothness of the quasi-interpolant; the corresponding error scales quadratically with the initial grid spacing. We validate our framework for tensor networks in one-, two-, and three-dimensional examples, including functions, derivatives, airfoil masks, and synthetic random fields such as three-dimensional turbulence. In particular, fractal fields can be generated directly in TT format with logarithmic memory and runtime. These results open a practical route to multiscale solvers, generative models, and geometry-aware algorithms on tensor-network and quantum platforms, with potential applications in scientific simulation, imaging, and real-time graphics.
\end{abstract}


\flushbottom
\maketitle

\section{Introduction}
\label{sec: introduction}

Upsampling, namely the construction of a finer-resolution representation from coarse data, is a fundamental operation in approximation theory, signal and image processing, and computer graphics. It includes interpolation, where the refined representation matches the given samples exactly, and quasi-interpolation, where it is reconstructed from local combinations of samples without requiring exact pointwise agreement, as well as more general forms of resampling and synthesis of fine-scale structure. These ideas underlie the approximation of multivariate functions \cite{multiva}, signal and image processing \cite{Schafer1973, Wang2022}, and procedural modeling in computer graphics \cite{Lagae2010,Perlin1985}. In practice, they appear in numerical PDE solvers \cite{Fasshauer1999, Hagstrom2015}, image resampling \cite{Parker1983,Keys1981}, and the construction of textures and noise fields, often termed \emph{synthetic noise} \cite{Lagae2010,Fournier1982}. A prominent extension is \textit{fractional Brownian motion} (fBm), or \textit{fractal noise} \cite{MidpointDisplacement,DiamondSquare,Perlin2002,Cook2005}, which has found applications in network traffic \cite{Leland1994}, hydrology \cite{Molz1997}, geophysical phenomena \cite{Mandelbrot1982}, porous media \cite{Jilesen2012}, turbulent dispersion \cite{Lilly2017}, and turbulence super-resolution \cite{Kim2008}. However, although classical upsampling schemes are highly effective in one dimension, in higher dimensions their cost typically grows exponentially with the number of dimensions. Even moderately refined meshes can therefore become intractable. This is another manifestation of the \textit{curse of dimensionality}, which affects a wide range of numerical problems, from machine learning and data science to the simulation of many-body quantum systems.

Two related representations exist that can alleviate the curse of dimensionality for certain structured data bases. 
The first one is given by tensor trains (TTs)~\cite{oseledets_tensor-train_2011}, also known as Matrix Product States in the quantum physics literature \cite{white,vidal}. This
provides a powerful low-rank factorization for multivariate arrays that can
effectively mitigate 
the curse of dimensionality in many practical cases, including the compression of structured functions~\cite{oseledets_constructive_2013,constructive_iteration,tindal2024}, PDE solvers~\cite{ttpde1,ttpde2,Peddinti2024,Markeeva2018-nv,Ali2025-aw},
turbulence modeling~\cite{Gourianov,stefano}, and multivariate analysis~\cite{GarcaRipoll2021}.
 The second one is given by quantum states, in the context of quantum computation. These can naturally  store an exponential amount of information.  Examples include states whose amplitudes encode probability distributions (\emph{q-samples}). Such q-samples underpin quantum applications boosting Monte Carlo and Markov-chain Monte Carlo \cite{Montanaro2015,Somma2008,Magniez2011,Paparo2014,Rebentrost2018} giving a quadratic speed up over classical Monte Carlo methods.
Interestingly, any quantum state admits a TT representation \cite{white,vidal}, where the ranks of the tensors depend directly on the entanglement of the state. Conversely, any (normalized) TT can be realized as a quantum state, prepared by a quantum circuit whose depth depend directly on the maximal rank of the tensors in the TT~\cite{seq_generation,Ran2020,Malz2024}.

On the tensor-network side, several strategies have been proposed to build TT representations of function-related tensors \cite{oseledets_tensor-train_2011, lindsey_multiscale_2024,qttchevishev,lubasch_multigrid_2018,ali_piecewise_2024}, among them sampling-based Tensor Cross Interpolation (\texttt{TT-Cross})~\cite{oseledets_tt-cross_2010,Savostyanov2011} has emerged as the state-of-the-art practical method for function approximation~\cite{NezFernndez2025,tree_cross}. However, it can overestimate TT ranks, its number of black-box evaluations may be the same order as the full tensor size, and its overall complexity still scales at least linearly with the number of TT-cores. In the quantum setting, a classical probability distribution over bit strings can be represented by a quantum state whose measurement outcomes reproduce the same distribution. This viewpoint underlies the quantum-sampling formulation of Aharonov and Ta-Shma \cite{qsamp_1}, the state-preparation method of Grover and Rudolph for efficiently integrable distributions \cite{grover_load}, and the construction by Low \emph{et al.} \cite{Low2014} of quantum states encoding Bayesian-network joint probability distributions.
Beyond probability distributions uploading on a quantum computer, only a few proposals address \emph{upsampling} of probability distributions \cite{chebtrans,MartnezdeLejarza2025,RamosCalderer2022}. However, these approaches either use global fits or patch-wise interpolation without boundary handling, leading to Gibbs-type oscillations near non-smooth features and limited control over upsampled smoothness.

We propose a unified upsampling framework with parallel realizations in \textit{tensor networks} (TN) and \textit{quantum states}. In the TN setting, it yields interpolation or quasi-interpolation of multivariate function values encoded on a coarse uniform grid, while in the quantum setting it yields quasi-interpolation of positive functions, such as probability distributions, encoded in a quantum state. In both cases, the coarse representation is refined to arbitrarily fine grids with controllable smoothness while preserving its tensor-network or quantum-state structure. The framework is agnostic to the underlying polynomial interpolation scheme, although for concreteness we focus on \emph{kernel polynomial interpolation} \cite{Keys1981,Thevenaz2000,unser}, which provides direct control over differentiability, accommodates nonperiodic boundaries in the TN setting, and yields derivatives at essentially no additional cost. In the TN setting, our method encodes one-dimensional functions in TT form in constant time and at fixed error, whereas the runtime of \texttt{TT-Cross} grows at least linearly with the number of cores; at 28 cores, we obtain roughly a three-order-of-magnitude speedup together with lower error. For a three-dimensional airfoil on grids with more than \(\sim 10^{9}\) points, \texttt{TT-Cross} fails to converge within a reasonable number of sweeps, while our method returns the encoded airfoil mask in constant time and at fixed error. The same framework also enables the generation of synthetic noise fields in TT format with logarithmic complexity and memory in the number of grid points; in particular, we construct a compressed three-dimensional synthetic turbulence field that reproduces Kolmogorov scaling and exhibits nontrivial intermittency in two different TN architectures. Further tensor-network applications, including image upsampling and the generation of one-dimensional noise functions and two-dimensional terrains, are presented in the Supplemental Material. In the quantum setting, we obtain a global quadratic quasi-interpolation error bound under periodic boundary conditions, and the construction extends to nonperiodic functions at the price of reduced boundary accuracy; the circuit depth scales logarithmically with the kernel degree and polynomially with the number of added qubits.

The paper is organized as follows. Section~\ref{sec: preliminaries} introduces the necessary background on tensor trains and interpolation. Sections~\ref{sec: QTT interpolation} and \ref{Sec:QC_SR} present our tensor-network and quantum state upsampling framework. Sections~\ref{sec: function_encoding}, \ref{sec: 3d_mask}, and \ref{Sec: synthetic_noise} show the main numerical results. Additional technical details are provided in Appendices~\ref{app:TTrep}--\ref{app: interpolation}. Appendix~\ref{app: synthetic noise} collects synthetic-noise algorithms in TT format. Further examples and metrics are given in the Supplemental Material.

\section{Preliminaries: Tensor Representations and Kernel Interpolation}
\label{sec: preliminaries}

\begin{figure*}[t!]
    \centering
    \includegraphics[width=1\textwidth]{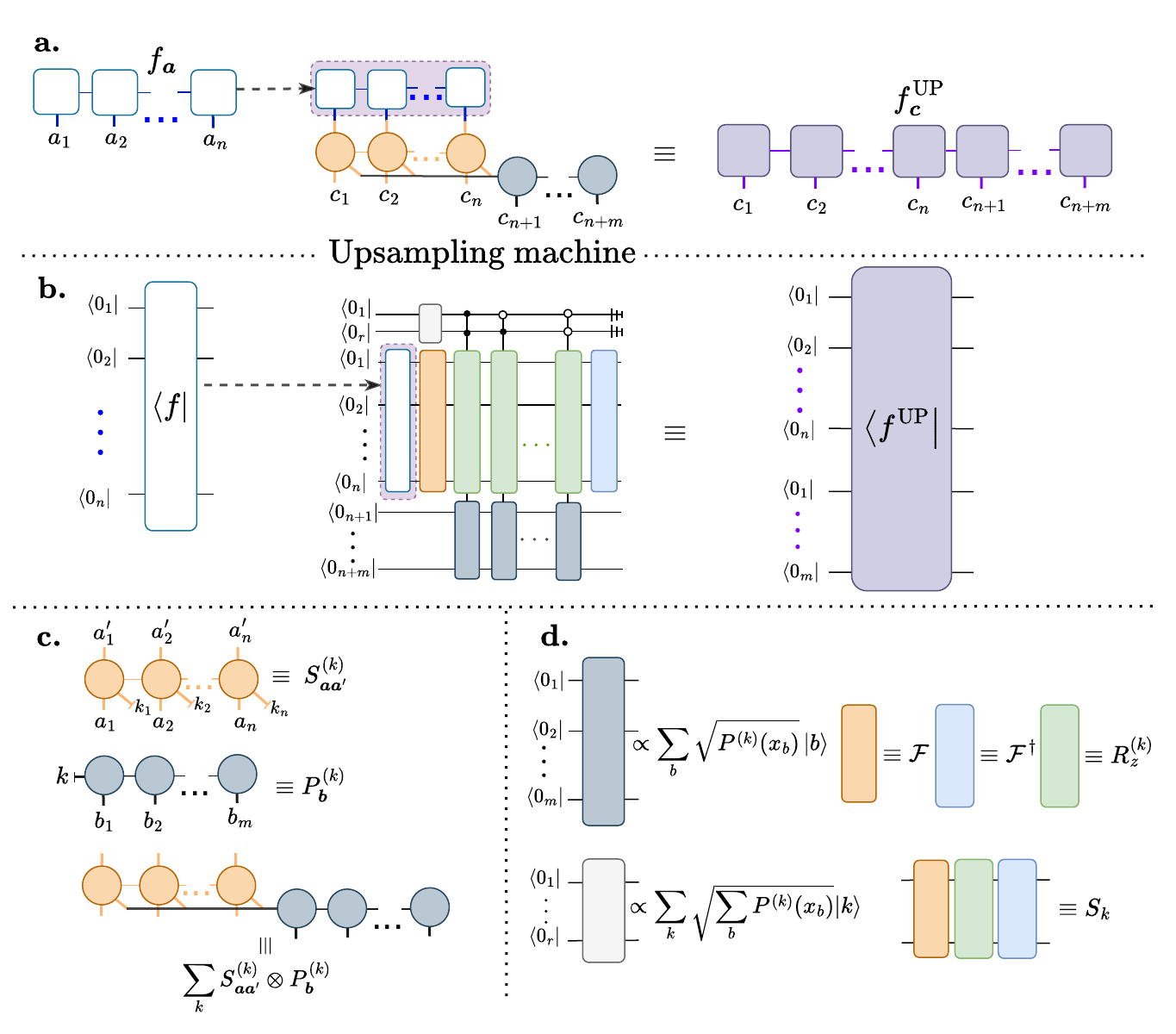}
\caption{
\textbf{Upsampling framework.} Schematic of the proposed multiscale upsampling framework for tensor-network and quantum-state representations. \textbf{(a)} Upsampling of an MPS, $f_{\boldsymbol a}$, encoding function values sampled on a dyadic grid via a mixed MPO--MPS tensor-train interpolation operator (TTI-O), which lifts an \(n\)-core MPS to an \(n+m\)-core representation with tunable smoothness of the upsampled function. \textbf{(b)} Upsampling (quasi-interpolation) of positive function values encoded in an \(n\)-qubit amplitude-encoded state onto \(m\) additional scales by a quantum circuit using \(r\) ancilla qubits. \textbf{(c)} Components of the TTI-O: shift matrices, \(S^{(k)}_{aa'}f_{a'} = f_{a+k}\), encoded as MPOs, and shifted polynomial kernels, \(P^{(k)}(x)=\phi\!\left(\frac{x}{h}-k\right)\), encoded as MPS as $P^{(k)}_{\boldsymbol{b}} \equiv P^{(k)}(x_{\boldsymbol{b}})$. Cut legs mean evaluating the tensor at that index. The indices $k$, $a$ and $a'$ are integers, bold indices are binary bit-strings, and indices with sub-indices take binary values. The full TTI-O is constructed as a superposition of tensor products of shift matrices and polynomial kernels. \textbf{(d)} Components of the quantum upsampler: \(\mathcal{F}\) denotes the quantum Fourier transform; \(R_{z}^{(k)}\) denotes a single layer of one-qubit \(z\)-rotation gates; in this case $\phi\!\left(\frac{x}{h}-k\right)$ is a positive kernel (e.g. a B-spline kernel) ; \(S_k\) is the unitary satisfying \(S_k \ket{a} = \ket{a+k \bmod 2^n}\); the \(m\)-qubit state-preparation gate is obtained through an MPS-to-QC encoder; and the \(r\)-qubit ancilla is prepared using a Hamming-weight encoder. }
    \label{fig:TTI}
\end{figure*}

We briefly summarize the tensor representations and interpolation tools that will be used throughout this work.  For a deeper explanation, review App. \ref{app:TTrep} and App. \ref{app: interpolation}.

A $d$--way tensor $\mathcal A \in \mathbb R^{n_{1}\times \cdots \times n_{d}}$ can be stored in \emph{Tensor Train} (TT) \cite{Oseledets2011} form as:
\begin{equation}
\label{eq: tt-from}
    \mathcal A_{\mathbf i}
    =
    \mathcal A_{i_{1}i_2\dots i_{d}}
    \;=\; 
G_{1}(i_{1})\,G_{2}(i_{2})\cdots G_{d}(i_{d}),
\end{equation}
where $\mathbf i=(i_1,\ldots,i_d)$ is a multi-index, with
$i_k\in\{1,\ldots,n_k\}$, and $G_{k}(i_{k}) \in \mathbb R^{r_{k-1}\times r_{k}}$ are matrix slices of three-dimensional tensors called \emph{cores} with $i_k \in \{1,\dots,n_k\}$ and \(r_0=r_d=1\). The matrix dimensions $r_k$ are called \emph{TT-ranks}. When all the physical indices, $i_k$, have size 2 the previous decomposition is called \textit{Quantics Tensor Train} (QTT) \cite{quantics}. 

It is possible to extend TT factorization to multidimensional linear operators,
called \textit{tensor train matrices} (TTM) \cite{Oseledets2010}, known as
Matrix Product Operators (MPO) in the physics literature \cite{earlympo,Pirvu2010},
as
\begin{align}
    \label{eq: tto}
    \mathcal O_{\mathbf i;\mathbf j}
    &=
    \mathcal O_{\,i_{1}\ldots i_{d};\,j_{1}\ldots j_{d}} \nonumber\\
&=\;
G_{1}(i_{1},j_{1})\,G_{2}(i_{2},j_{2})\cdots G_{d}(i_{d},j_{d}),
\end{align}
where $\mathbf i=(i_1,\ldots,i_d)$ and $\mathbf j=(j_1,\ldots,j_d)$ are the
row and column multi-indices, respectively. Here $i_k\in\{1,\dots,n_k\}$,
$j_k\in\{1,\dots,m_k\}$, and each
$G_k(i_k,j_k)\in\mathbb{R}^{r_{k-1}\times r_k}$ is a matrix slice of a
four-dimensional tensor, with $r_0=r_d=1$.

\begin{figure*}[t!]
    \centering
    \includegraphics[width=\textwidth]{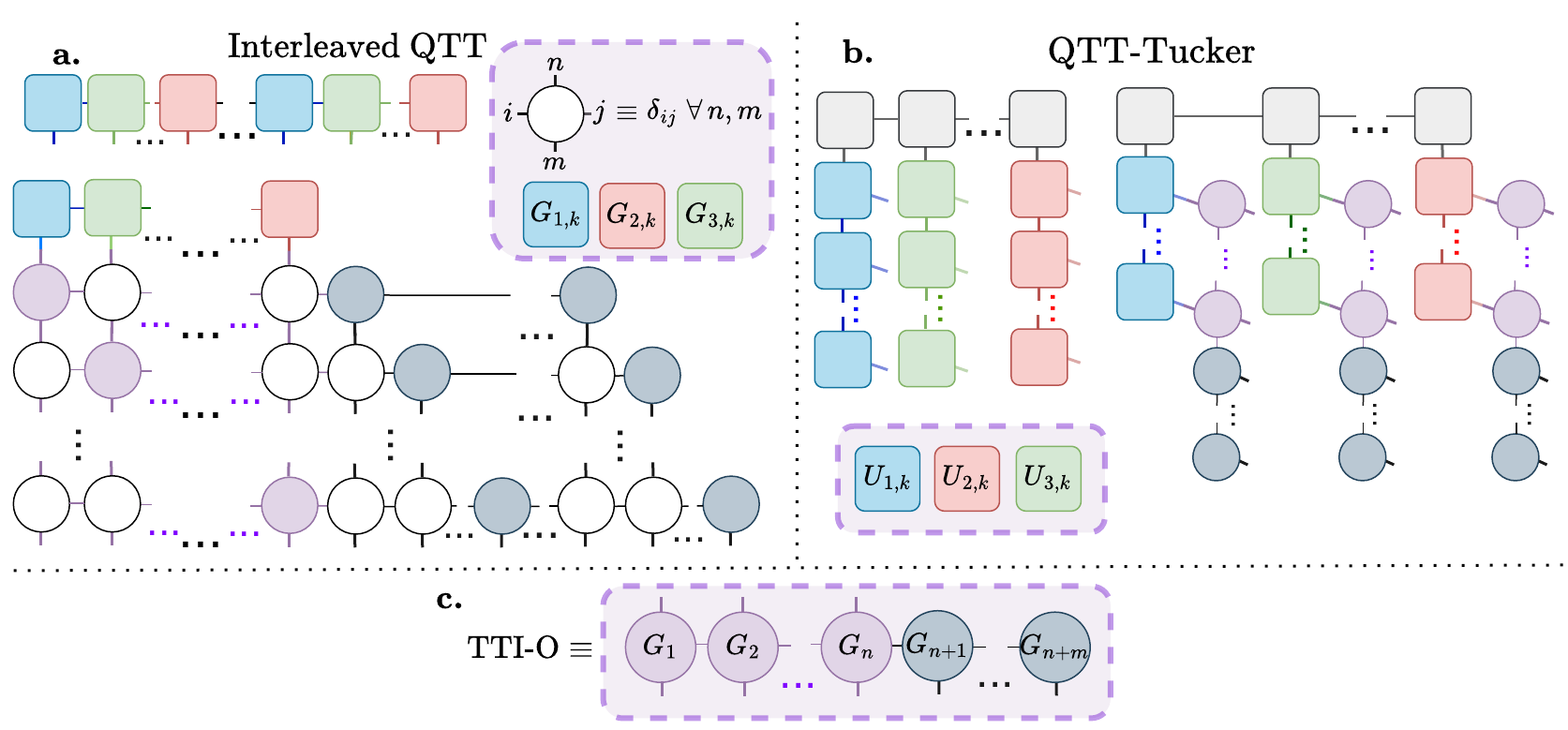}
\caption{\textbf{Multi-dimensional Tensor Train Interpolation (TTI).} \(\mathbf{(a)}\) \(d\)-dimensional interpolation in the \textit{QTT-I} format: each TTI-O acts along one dimension, with identity padding on the remaining dimensions, yielding a product of \(d\) coupled one-dimensional TTI-O operators. \(\mathbf{(b)}\) \(d\)-dimensional interpolation in the \textit{QTT-T} format: each TTI-O acts independently on its corresponding dimension. Different colors denote different dimensions, and \(k\) labels the scales. Here, \(G_{i,k}\) denotes the cores in the \textit{QTT-I} representation, while \(U_{i,k}\) denotes the unitary cores in the \textit{QTT-T} representation. \textbf{(c)} Graphical representation of the TTI-O defined in Eq.~\ref{eq: TTIO}.}
    \label{fig:ndtt}
\end{figure*}

Several QTT-like extensions exist to encode multi-dimensional tensors, $\mathcal A \in \mathbb R^{n_{1}\times \cdots \times n_{d}}$, where $n_k = 2^{\alpha_k}$. For simplicity, let us assume $\alpha_k = \alpha$ for all $k$. In this work, we consider two such formats: \textit{QTT-interleaved} (QTT-I) \cite{Ye2024} and \textit{QTT-Tucker} (QTT-T) \cite{Tucker1966,DeLathauwer2000}. In both cases, each dimension, $n_k$, is decomposed as a binary multi-index through its binary expansion $n_k  = \sum_{i=1}^m a_{k,i} 2^{i}$, thus $n_k \rightarrow a_{k,1}\dots a_{k,m}$. We will refer to each binary power $i$ as a \textit{scale}. The QTT-I format orders the TT physical indices scale by scale. 
For each scale $i$, the binary indices $a_{k,i}$ associated with the different physical dimensions are grouped in ascending order, which results in the multi-index ordering $
a_{1,1}\ldots a_{d,1}
a_{1,2}\ldots a_{d,2}
\ldots
a_{1,m}\ldots a_{d,m}
$. In contrast, the QTT-T format uses a nested representation: a global TT decomposition  that separates the dimensions $n_k$, while each dimension is further decomposed into a QTT decomposition. A schematic overview of these two encodings is shown in Fig.~\ref{fig:ndtt}\textbf{a} and Fig.~\ref{fig:ndtt}\textbf{b}, where different colors represent different dimensions. Further technical details are provided in Appendix~\ref{app:TTrep}.

Different schemes for constructing a TT representation of a tensor include: hierarchical SVD-based methods (\texttt{TT-SVD})~\cite{oseledets_tensor-train_2011}, multiscale interpolative QTT schemes~\cite{lindsey_multiscale_2024}, Chebyshev-based approaches~\cite{qttchevishev}, MPO-based multigrid refinements combined with \texttt{DMRG}-like optimization~\cite{lubasch_multigrid_2018}, QTT constructions with embedded piecewise polynomial bases~\cite{ali_piecewise_2024}, and the aforementioned Tensor Cross Interpolation (\texttt{TT-Cross})~\cite{oseledets_tt-cross_2010,Savostyanov2011}. However, \texttt{TT-SVD} still requires access to the full tensor and is therefore limited by available RAM; multiscale interpolative schemes typically require function evaluations on non-regular grids; and \texttt{TT-Cross} may fail to recover an accurate approximation within a reasonable amount of time.

Turning to interpolation, we use \textit{kernel convolution interpolation} \cite{Keys1981,Thevenaz2000,unser} with compactly supported polynomial kernels, whose degree and support control both the smoothness of the interpolant and the interpolation error. For simplicity, we will assume periodic boundary handling. Let $f:[0,1)\rightarrow\mathbb{R}$ take values on a uniform coarse grid of $N$ points, $\{x_a\}_{a=0}^{N-1}$, with step size $h=x_{a+1}-x_a=1/N$, and let $\phi:\mathbb{R}\to\mathbb{R}$ be an interpolation kernel such that $\phi(x) = 0$ for $|x|>\frac{q+1}{2}$, where $q$ is a natural number. Then, on the interval $[x_a,x_{a+1})$, the interpolant is given by
\begin{equation}
F(x)=\sum_{k=-\lceil q/2\rceil}^{\lfloor q/2\rfloor} f(x_{a+k})\,\phi\!\left(\frac{x}{h}-k\right).
\end{equation}
Notice that interpolation requires $F(x_a)=f(x_a)$. More generally, it is possible to construct kernels that give approximants with $F(x_a)\neq f(x_a)$ and $\norm{F - f}_2 \leq C h^m$, where $m$ and $C$ depend on the regularity of $f$ and the degree of the kernel; this is known as \emph{quasi-interpolation} \cite{Boor1990,Sablonniere2005}.  Furthermore, high dimensional kernels factorize as a product of one-dimensional kernels, $\boldsymbol{\phi}(\mathbf{x})=\prod_{j=1}^d \phi(x_j)$, so multidimensional interpolation is realized as a sequence of one-dimensional interpolations.

\section{Tensor Train Interpolation}
\label{sec: QTT interpolation}

In this section, we present our tensor-network construction for interpolation. To begin with, we will explain the technique in 1D and the generalization to the multi-dimensional case follows immediately. Let $\phi$ be a polynomial kernel with finite support. For simplicity, we assume periodic boundary conditions. Consider a function $f:[0,1)\rightarrow\mathbb{R}$ sampled on a uniform grid of $2^n$ points, $\{x_a\}_{a=0}^{2^n-1}$, and suppose that we wish to refine its representation to a finer uniform grid of $2^{n+m}$ points. We write the index $a$ in binary form as $a=\sum_{i=1}^{n}2^{\,n-i}a_i$, with $a_i\in\{0,1\}$, and denote the associated multi-index by $\boldsymbol{a}=a_1a_2\dots a_n$. The corresponding grid points are labeled as $x_a=2^{-n}a$. We denote the QTT representation on $n$ cores by $f_{\boldsymbol a}\equiv f(x_{\boldsymbol a})$. The additional refined scales are labeled by a second multi-index $\boldsymbol b=b_1b_2\dots b_m$, with associated coordinate $x_b=2^{-m}b$.

Since the kernel $\phi$ has finite support, let $q$ be the number of neighbors that contribute to the interpolation. Therefore, we can write the interpolation as a superposition of $q+1$ polynomials of degree $p$ weighted by the corresponding shifted function value, resulting in:

\begin{equation}
\label{eq: ttint}
\begin{split}
   \text{TTI}(f_{\boldsymbol{a}}, \mathbf{\phi}_{\boldsymbol{b}}) &\equiv \sum_{\boldsymbol{a'}} \left(  \sum_{k}  S^{(k)}_{\boldsymbol{a}\boldsymbol{a}'} \otimes P^{(k)}_{\boldsymbol{b}}  \right)   f_{\boldsymbol{a}'}\\ 
\end{split}
\end{equation}

\noindent where $k \in {-\lfloor q/2\rfloor,...,\lfloor (q-1)/2\rfloor}$, $S^{(k)}_{\boldsymbol{a}\boldsymbol{a}'}$ are shift matrices, $S^{(k)}_{{a}{a}'}  f_{{a}'} = f_{a+k}$, and $P^{(k)}(x) = \phi(x/h-k)$ are polynomials of degree $p$ defined on $[0,1)$, encoded as QTTs as $P^{(k)}_{\boldsymbol{b}} \equiv P^{(k)}(x_{\boldsymbol{b}})$. Shift matrices are rank 2 MPOs \cite{Kazeev2013} and encoding polynomials in QTT format is a known construction \cite{oseledets_constructive_2013}. Eq.~\eqref{eq: ttint} interpolates $f$ on $m$ new sub-scales.  
We can gather the shift matrices with their corresponding polynomial as an MPO with ranks bounded by $q+1$ on the first $n$ legs and an MPS on the last $m$ interpolated legs with ranks bounded by $p+1$, see Fig.~\ref{fig:TTI}. 
This MPO-MPS operator performs the interpolation over the finer grid.

Let \(M=n+m\) be the total number of scales, the first \(n\) scales carry operator legs while the remaining \(m\) scales are vector legs. Let's define $ c_k =(a_k,\, a'_k)$ if $ k \le n$ and $ c_k = b_{k-n}$ if $ k>n$, where $a_k$, $a'_k$ and $b_k$ take binary values. Therefore, we can write the 1D TTI operator (TTI-O) as:

\begin{equation}
\label{eq: TTIO}
    \text{TTI-O} = G_1( c_1) \cdots G_M( c_M)
\end{equation}
with $G_k(c_k) \in \mathbb{R}^{r_{k-1} \times r_k}$ and $r_0=r_M=1$. This operator acts on an MPS as a normal MPO-MPS contraction over the first $n$ legs, see Fig.~\ref{fig:TTI}\textbf{a}.

Multi‐dimensional convolution interpolation is realized as an iterated sequence of 1D interpolations, refining one coordinate at a time. For QTT-T this is trivial, since we can apply TTI to each dimension independently, see Fig.~\ref{fig:ndtt}\textbf{b}. On the other hand, for QTT-I we extend the 1D TTI–O to a \(d\)-dimensional operator. This construction was done in \cite{Markeeva2020} for two-dimensional operators, but the generalization to $d$-D is straightforward. We build the cores of the multidimensional TTI-O using the 1D cores and padding identities such that the cores labeled by the same dimension $m$ are acted by the same TTI-O, see Fig.~\ref{fig:ndtt}.\textbf{a}. Explicitly, a core labeled by dimension $m \in \{1,\ldots,d\}$ and scale $k\in \{1,\ldots,M\}$ takes the form: 

\begin{equation}
\label{eq:multicore}
    G_{m,k}( c_{m,k}) = \mathds{1}^{\otimes(m-1)}_{r_{k-1}}\otimes G_k( c_{m,k}) \otimes \mathds{1}^{\otimes(d-m)}_{r_k}.
\end{equation}

\noindent Here, $\mathds{1}_r$ denotes the $r \times r$ identity matrix.  With this we can immediately see that the QTT-T representation gives a better compression since the rank tails of each QTT-T leg are always bounded by $p+1$, while for QTT-I the rank tails are bounded by $(p+1)^d$, since the multidimensional TTI-O can be seen as the product of $d$ one dimensional TTI-O.

Moreover, our framework is not restricted to kernel-based polynomial interpolation, other interpolation and quasi-interpolation methods, such as Lagrange interpolation \cite{stoer_introduction_2002}, can be incorporated as well; see Supplemental Material Sec.~IIA for \(C^0\) interpolants and Sec.~II for a broader review. The construction also extends straightforwardly to non-periodic functions by modifying the shift matrices to include the appropriate boundary terms. In the tensor-network setting, these corrections are implemented through element-wise matrix additions or subtractions, each representable as a rank-1 MPO.

\section{Shallow Quantum Upsampling}
\label{Sec:QC_SR}

\begin{figure*}[t!]
    \centering
    \includegraphics[width=\textwidth]{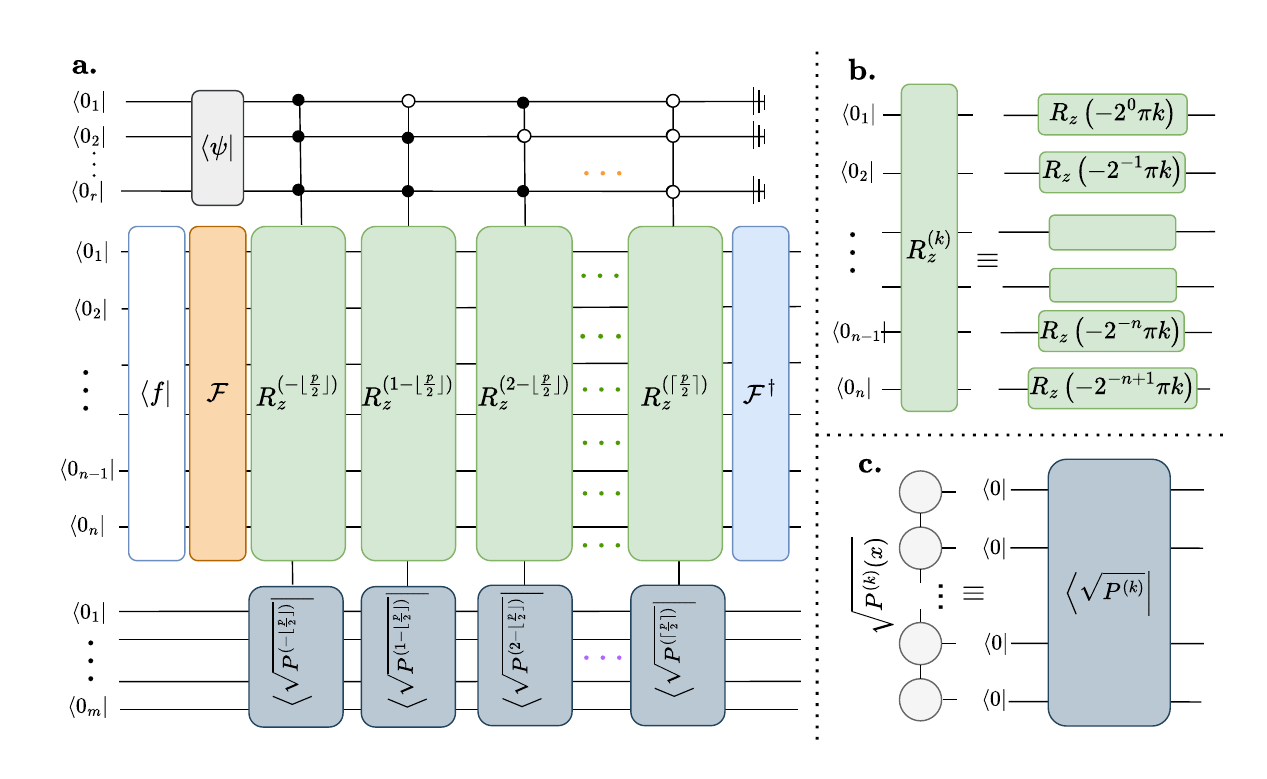}
\caption{\textbf{Quantum-Circuit Spline Quasi-Interpolation.}
A shallow quantum routine to \emph{refine} a function encoded in a quantum state.
\textbf{(a)} Full quasi-interpolation circuit.
The circuit first prepares the ancilla state $\ket{\psi}$ (see Eq.~\eqref{eq:psi}). After applying the QFT $\mathcal{F}$, the $\ket{k}$ register controls the shift operators (see Eq.~\eqref{eq:qc-shift-qft}) and the associated polynomial blocks. Black dots denote controls and white dots denote anti-controls. Grounded wires indicate indicate a partial trace over those qubits. In the Fourier domain, the shift operators are implemented as $R_z^{(k)}$ rotations, and the circuit prepares the polynomial states $\ket{\sqrt{P^{(k)}}}$.
\textbf{(b)} Decomposition of the diagonal gate $R_z^{(k)}$ into single-qubit rotations.
\textbf{(c)} Schematic of the MPS-to-quantum-circuit (MPS--QC) encoder (see~\cite{Ran2020,Malz2024}).}
    \label{fig:mps-qc}
\end{figure*}

In this section, we show how to upsample a positive function, such as a probability density, encoded in the probabilities of a quantum state, using the convolution-kernel framework of Sec.~\ref{sec: preliminaries}. Because upsampling occurs in the probabilities, avoiding post-processing of the quantum samples requires a positive kernel. The construction is therefore restricted to positive quasi-interpolation kernels, which yield a quasi-interpolant with quadratic approximation error in the grid spacing. A schematic of this method is shown in Fig.~\ref{fig:mps-qc}.

Let $f:[0,1)\to\mathbb{R}^+$ be a positive $C^1$ function encoded in an $n$-qubit quantum state as:

\begin{equation}
\label{eq:qc-fstate}
    \ket{f}
    =
    \frac{1}{\sqrt{Z_n}}
    \sum_{a=0}^{2^n-1} \sqrt{f_a}\,\ket{a},
    \qquad
    Z_n := \sum_{a=0}^{2^n-1} f_a .
\end{equation}

\noindent where $a=\sum_{j=1}^{n}2^{n-j}a_j$, $x_a = 2^{-n}a$ and $f_a = f(x_a)$. As before, we impose periodic boundary conditions, so the values outside \([0,1)\) are wrapped around according to $f(1+x)=f(x)$. This includes, for example, periodic functions and symmetric extensions. Moreover, the $n$-qubit register is labeled by the multi-index $\mathbf{a}=a_1 \ldots a_n$ and to refine $\ket{f}$ by $m$ additional scales, we introduce an $m$-qubit register $\mathbf{b}=b_1\cdots b_m$, with $b=\sum_{\ell=1}^{m}2^{m-\ell}b_\ell$ and local coordinate $x_b=2^{-m}b$. The refined basis index is then $\mathbf{e}=(\mathbf{a},\mathbf{b})$, or equivalently $e=2^m a+b$.

To construct the upsampled state $\ket{\tilde{f}}$, we use a degree-$p$ B-spline centered kernel $\beta^{(p)}(x)$ and define the shifted polynomial pieces $P^{(k)}(x)=\beta^{(p)}(x-k)$, $k=-\lfloor p/2\rfloor,\dots,\lceil p/2\rceil$.
The coefficients of $\ket{\tilde{f}}$ over the refined basis $\{\ket{e}\}_{c=0}^{2^{n+m}-1}$ should satisfy:

\begin{equation}
\label{eq:qc_int}
    \braket{e}{\tilde{f}} \propto \left( \sum_{k = -\lfloor p/2\rfloor}^{\lceil  p/2\rceil} f_{a+k} P^{(k)}(x_b) \right)^{1/2}
\end{equation}
In order to achieve this, we construct a quantum gate $U_I$ such that it acts on $q=\lceil\log_2(p+1)\rceil$ ancilla qubits, the input $n$-qubit state $\ket{f}$,  and the new $m$-qubit register as:

\begin{equation}
\label{eq:qc-map}
    U_I \left( \ket{0}_q\,\ket{f}\,\ket{0}_m \right)
    =
    \sum_{k = -\lfloor p/2\rfloor}^{\lceil  p/2\rceil}
    \alpha_{k}\, \ket{k}\,
    S_{k}^\dagger\left( \ket{f} \right)\,
    \ket{\sqrt{P^{(k)}}}.
\end{equation}

\noindent Here $\alpha_k$ are proportional to the square root of the $\ell_1$ norms of the polynomials $P^{(k)}(x)$, $\alpha_{k}:=\sqrt{2^{-m}\sum_{b} P^{(k)}(x_b)} $. Ancilla states $\ket{k}$ with negative $k$ are defined as $\ket{k \mod p+1}$ . The operators $S_k\ket{a}=\ket{a+k \bmod 2^n}$ are modular shifts and the states $\ket{\sqrt{P^{(k)}}}$ are given by

\begin{equation}
\label{eq:qc-kernelstate}
 \ket{\sqrt{P^{(k)}}}=   \frac{1}{\alpha_k}    \sum_{b}\sqrt{2^{-m}\,P^{(k)}(x_b)}\,\ket{b} .
\end{equation}

Moreover, since $S_k$ is diagonal in the Fourier basis~\cite{Shakeel2020-yv}, it admits the decomposition
\begin{equation}
\label{eq:qc-shift-qft}
    S_k
    =
    \mathcal{F}^{\dagger}
    \left(\bigotimes_{j=1}^{n} R_z\!\left(\frac{\pi k}{2^{j-1}}\right)\right)
    \mathcal{F},
\end{equation}
where $\mathcal{F}$ is the quantum Fourier transform (QFT). Therefore, shifts reduce to products of one-qubit rotations in Fourier space.


In the following, we describe the components of $U_I$, as illustrated in Fig.~\ref{fig:mps-qc}\textbf{a}.
First, we prepare on the $q$ ancilla qubits the state
\begin{equation}
\label{eq:psi}
    \ket{\psi}=\sum_{k = -\lfloor p/2\rfloor}^{\lceil  p/2\rceil} \alpha_k\ket{k}.
\end{equation}
Since the B-spline kernel forms a partition of unity, $\ket{\psi}$ is normalized.
Second, to implement the controlled modular shifts $S_k^\dagger$ on $\ket{f}$ with $\ket{k}$ as controls,
we apply  $\mathcal{F}$, so that each controlled shift reduces to controlled single-qubit $R_z$
rotations in Fourier space.
Third, using the same controls, the circuit prepares the $m$-qubit state $\ket{\sqrt{P^{(k)}}}$.
Finally, we apply $\mathcal{F}^\dagger$ to return to the computational basis.
The control state $\ket{k}$ is implemented by matching the binary expansion of $k$: starting from $\ket{0}_m$, qubits corresponding to binary $1$ are controls, while qubits corresponding to binary $0$ are converted into controls by applying an $X$ gate before and after the controlled operation.
Measuring only the $(n+m)$-qubits, $U_I \left( \ket{0}_q\,\ket{f}\,\ket{0}_m \right)$ satisfies Eq.~\eqref{eq:qc_int}.

The ancilla state $\ket{\psi}$ can be prepared efficiently~\cite{Farias2025}, and controlled one- and two-qubit rotations can be optimized following~\cite{Vale2024}. The controlled polynomial states are compiled using an MPS--QC encoder~\cite{Ran2020,Malz2024}.
The dominant cost is typically the QFT; nevertheless, this is more efficient than implementing a controlled modular adder.
The differentiability of the upsampled function is tuned by the degree $p$ of the B-spline kernel.
In practice, the kernel states $\ket{\sqrt{P^{(k)}}}$ can be obtained from \texttt{TT-SVD}; \texttt{TT-Cross} is also viable for a moderate number of added scales, while \texttt{TTI} is preferable when very fine grids are required.
After TT rounding, we observe that the maximum TT-rank remains bounded by $10$ regardless of $p$; therefore, the circuit depth depends mainly on the number of ancillas $q$. Furthermore, due to the quasi-interpolation scheme used, the quasi-interpolation error $\mathcal{O}(2^{-2n})$ does not depend on the degree of the polynomial kernel $p$.

\begin{figure*}[t!]
    \centering
\includegraphics[width=\textwidth]{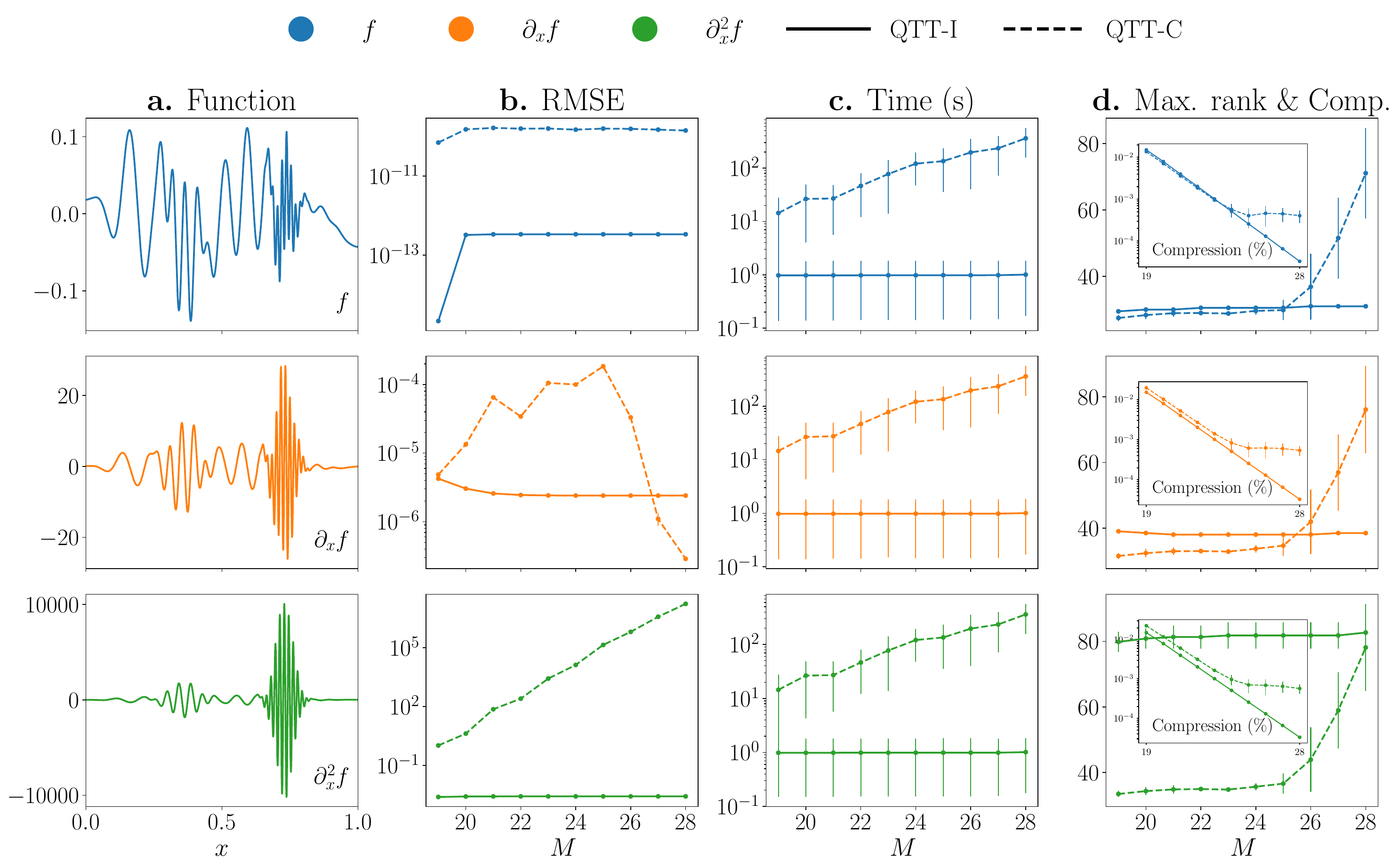}
\caption{\textbf{1D function encoding.}
Demonstration of tensor-train interpolation (TTI) in one dimension. 
A \(C^{2}\) function \(f:[0,1) \rightarrow \mathbb{R}\) and its derivatives are represented in QTT format. We denote the number of scales (cores) by $M$.
Starting from a coarse QTT representation of $f$ with \(18\) cores (grid spacing \(h = 2^{-18} \approx 10^{-6}\)), we refine the representation up to \(28\) cores using TTI; this interpolated representation is denoted \textit{QTT-I}. 
For comparison, the same function is also encoded with \texttt{TT-Cross}, denoted \textit{QTT-C}. 
\textbf{(a)} Function and derivatives. 
\textbf{(b)} Root-mean-square error (RMSE) of the function and its first two derivatives, evaluated by TT sampling~\cite{Ferris2012-mg}. 
Derivatives in \textit{QTT-I} are obtained analytically from the interpolation, whereas those in \textit{QTT-C} are computed with a finite difference MPO. 
\textbf{(c)} Runtime of TTI and \texttt{TT-Cross} for a fixed target precision. 
\textbf{(d)} Maximum QTT rank for both methods, with an inset showing the corresponding compression ratio. 
}

\label{fig:1d_function}
\end{figure*}

\section{Applications}
\label{sec: results}

\begin{figure*}[t!]
    \centering
\includegraphics[width=\textwidth]{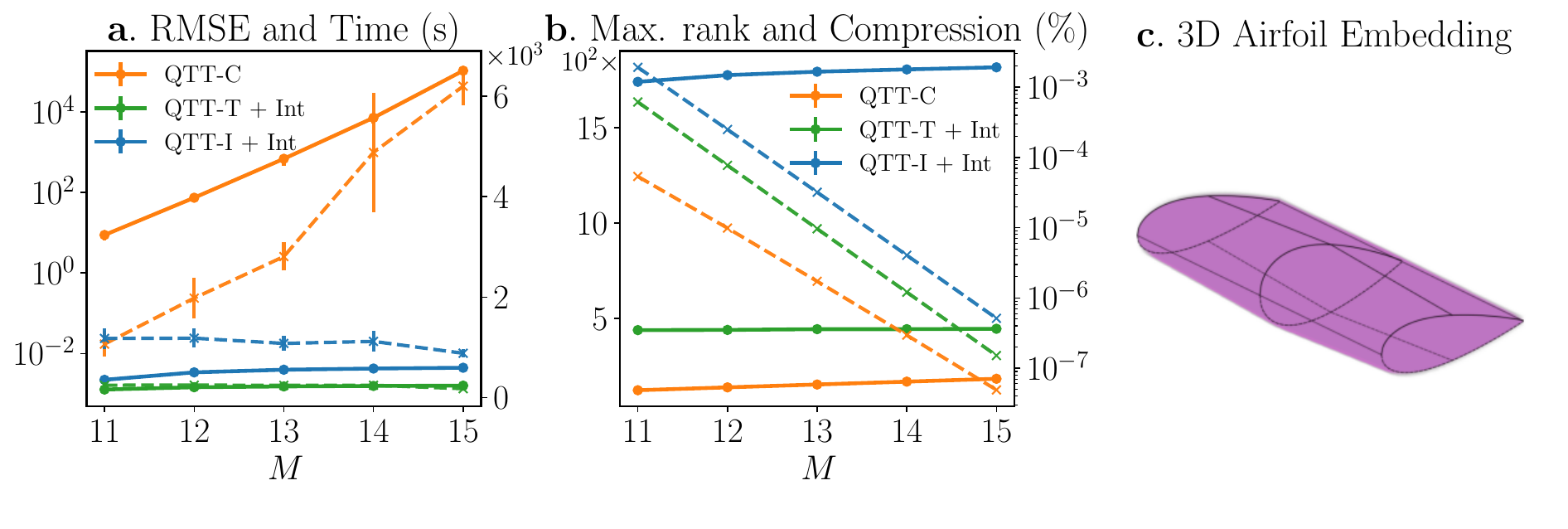}
\caption{\textbf{3D soft-mask encoding.}
Three-dimensional soft indicator function of a tampered airfoil embedded in a computational domain that is four times larger than the object. 
The results compare three tensor-network approaches: \texttt{TT-Cross} (\textit{QTT-C}), tensor-train interpolation in the QTT-interleaved architecture (\textit{QTT-I}), and tensor-train interpolation in the QTT–Tucker architecture (\textit{QTT-T}). The interpolation kernel is a $C^1$ cubic kernel. The number of scales per dimension (one third of the number of cores) is denoted by $M$.
\textbf{(a)} Runtime of the three methods (solid lines), together with the root-mean-square error (RMSE) computed via TT sampling~\cite{Ferris2012-mg} (dashed lines). 
\textbf{(b)}  Maximum TT rank (solid lines) and compression ratio (dashed lines) for each tensor-network architecture. 
\textbf{(c)} Visualization of the softened mask of the airfoil in the full 3D domain. 
Panels~\textbf{(a)} and~\textbf{(b)} demonstrate that our interpolation-based method yields constant runtime, maximum rank, and error independen of the number of upsampled scales. 
}

\label{fig:3d_wing}
\end{figure*}

We illustrate our upsamplig framework for tensor-networks on two representative tasks: \emph{function/mask encoding}, and \emph{synthetic noise}. We present a few examples on the main text, but several more can be found in the 
Supplemental Material Sec. III. For the first one, see Sec. \ref{sec: function_encoding}, we compare \texttt{TT-Cross} versus direct interpolation in 1D, 2D and 3D, reporting runtime, maximum bond dimension, compression and approximation error as functions of grid resolution.
For the second, see Sec. \ref{Sec: synthetic_noise}, we generate 1D, 2D and 3D noise fields by combining pseudorandom values/gradients with our upsampling method, yielding scale‐invariant textures with tunable smoothness. Detailed algorithms are given in Appendix~\ref{app: synthetic noise} and Appendix \ref{App:turbulence}.

\subsection{Soft Masks}
\label{sec: function_encoding}

Here we show how to \textit{boost} the QTT construction of soft indicator functions, \textit{soft masks}, over exponentially fine grids, i.e. high number of TT-cores. 
A mask is an indicator function of a region of interest, taking the value $1$ inside the domain and $0$ outside, while a \textit{soft} mask smooths this discontinuous jump (e.g., via a bump function) taking values between $0$ and $1$.
First, we build a coarse QTT representation of a function using \texttt{TT-SVD} \cite{oseledets_tensor-train_2011} (although \texttt{TT-Cross} \cite{oseledets_tt-cross_2010} could be used as well) as our base for interpolation and then use TTI to interpolate it into finer grids.

To illustrate TTI, we consider examples in one, two, and three dimensions. 
The 1D example demonstrates the effectiveness of TTI by upsampling a function with scale-dependent oscillatory behavior while also approximating its derivatives at no additional cost.
In two dimensions, we consider two simple but representative cases: a correlated Gaussian distribution on a domain containing $99.999\%$ of its mass, and a soft mask for an airfoil/circle centered in a rectangular domain, see Supplemental Material Sec.~III. 
For the 3D example, we encode a softened indicator function of a tampered airfoil and compare the performance of \texttt{TT-Cross} against our TTI procedure. 
All benchmarks compare the accuracy and compression of TTI with the implementation of TT-Cross in \texttt{Teneva}~\cite{chertkov2024teneva}. 
The tensor-train computations were performed using \texttt{torchTT}~\cite{torchtt}.

\subsubsection{1D function}

We begin with a $C^2$ function that exhibits different oscillatory behaviors across its domain. 
Because of these oscillations, accurate interpolation requires a sufficiently fine initial sampling, with at least two samples inside each oscillation period. 

We first construct a QTT representation with 18 cores, i.e. a coarse grid spacing of $h=2^{-18}$, using \texttt{TT-SVD}, and then refine it with TTI using a $C^1$ cubic kernel with $\mathcal{O}(h^3)$ interpolation error, see Supplemental Material Sec.~II. Therefore, expected interpolation error is \(\mathcal{O}(10^{-17}\delta)\) where $\delta$ depends on the second derivative of the function.
The first derivative is obtained by differentiating each polynomial piece $P^{(k)}(x)$ in Eq.~\ref{eq: ttint}. The interpolation error on the first derivative is $\mathcal{O}(10^{-11}\delta)$.
For the second derivative, we use a $C^2$ cubic B-spline kernel, see Supplemental Material Sec.~II, and differentiate the corresponding pieces twice. This kernel has a quasi-interpolation error of \(\mathcal{O}(10^{-6}\eta)\), where \(\eta\) depends on the fourth derivative of the function.
As a baseline, on each refined scale we construct a QTT representation using \texttt{TT-Cross} and approximate derivatives applying a finite-difference MPO \cite{Kazeev_low-rank_2012}.

As Fig.~\ref{fig:1d_function}\textbf{b} shows, once a function is encoded in QTT form at sufficiently high resolution, upsampling it to arbitrarily fine grids becomes straightforward with TTI.  The first row shows that the interpolation error for the function  is of order \(\mathcal{O}(10^{-13})\). For the first derivative, the error is \(\mathcal{O}(10^{-6})\), while for the second derivative it is $\mathcal{O}(10^{-2})$, as shown in the second and third rows, respectively. For a fixed number of sweeps, the \texttt{TT-Cross} error increases with the number of cores, indicating that it fails to find an accurate representation of the second derivative. In contrast, the TTI $\ell^2$ error is controlled by analytical interpolation bounds, so the approximation error for both the function and its derivatives is constant regardless of the number of upsampled scales. Moreover, Fig.~\ref{fig:1d_function}\textbf{c} shows that the runtime of \texttt{TT-Cross} grows at least linearly with the number of cores. In contrast, TTI runs in constant time plus minor corrections, since it only requires the initial QTT encoding and a final TT-rounding step on the coarse scales. Finally, Fig.~\ref{fig:1d_function}\textbf{d} shows that, for a large number of cores, \texttt{TT-Cross} tends to overestimate the TT-ranks, while TTI has a constant max rank resulting in an exponential compression.

\begin{figure*}[t!]
    \centering
\includegraphics[width=\textwidth]{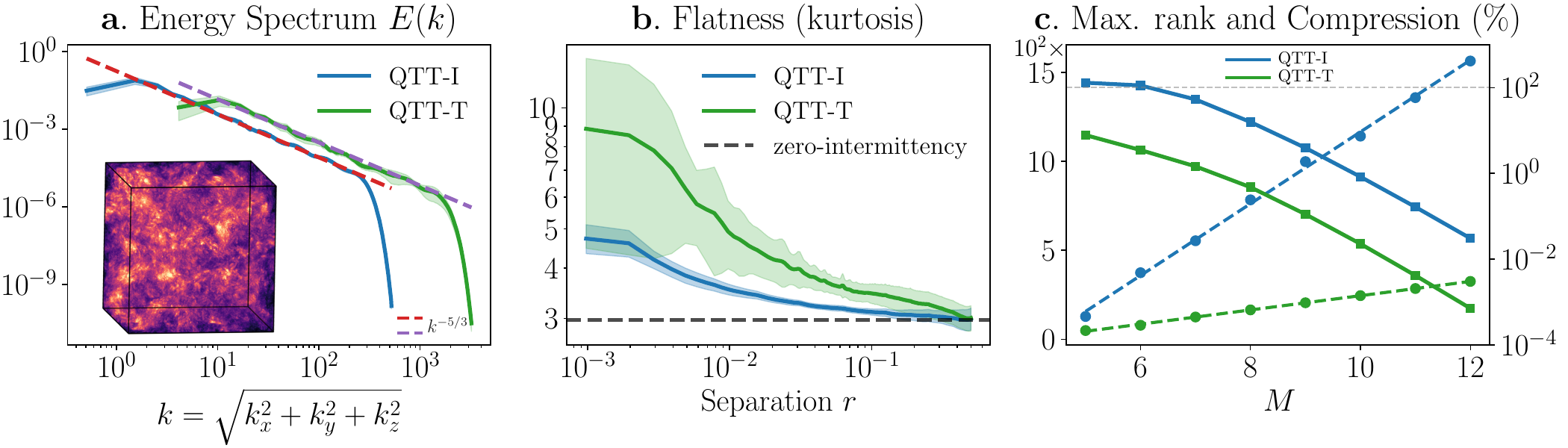}
\caption{ \textbf{Synthetic turbulence.} 3D turbulence metrics for two tensor-network variants: \emph{QTT-Interleaved} (QTT-I) and \emph{QTT-Tucker} (QTT-T).  \textbf{(a)} Energy spectrum $E(k)$ vs.\ wavenumber $k=|\mathbf{k}|$ (log–log); dashed red and purple reference lines (QTT-I and QTT-T, respectively) indicate $\propto k^{-5/3}$, confirming the expected Kolmogorov inertial-range scaling; inset: representative 3D snapshot of the velocity magnitude. \textbf{(b)} Flatness (kurtosis) of velocity increments as a function of separation $r$ (log–log); systematic departures from a constant baseline highlight intermittency (non-Gaussian fluctuations) across scales. \textbf{(c)} Maximum bond dimension $r_{\text{max}}$ versus the number of scales $M$ (circles, dashed lines; one third of the TN cores), showing approximately linear growth. Compression ratios (squares, solid lines; TN parameters divided by grid points) are plotted on the secondary axis to quantify storage efficiency. Results are averaged over 20 random seeds with a fixed number of scales $M=10$ per dimension; shaded regions indicate $\pm 1$ standard deviation.
}
\label{fig:synthetic_turbulence}
\end{figure*}

\subsubsection{3D  masks}
\label{sec: 3d_mask}

The idea of encoding indicator functions as QTTs was first introduced in \cite{Peddinti2024}, where they incorporated 2D objects into a complete \textit{quantum-inspired} pipeline for \textit{computational fluid dynamics}. In \cite{Peddinti2024}, they show that \texttt{TT-Cross} performs better on softened indicator functions, allowing for the correct encoding (no artifacts) of non-slip boundary conditions around the object's boundaries. Moreover, a related work on compressing 3D objects was presented in \cite{Boyko2020}, where the authors show that representing \textit{Truncated Signed Distance Functions} (TSDFs) in TT format preserves visual rendering quality.

Following the idea of softening an indicator function in order to have a low-rank object \cite{Peddinti2024}, we softened the indicator function of a tampered airplane wing, and then we embed it on a lattice four times larger (per dimension) than it, leaving sufficient room for vortical structures relevant to fluid simulation. As before, we start with a TT representation of the wing on a coarse grid with $2^{10}$ points per dimension. Then, we resolute new scales with a $C^1$ cubic kernel applied across all dimensions. We compare the running time of TTI on both QTT-interleaved and QTT-Tucker against \texttt{TT-Cross}.
As mentioned in Sec.~\ref{sec: QTT interpolation} and shown in Fig.~\ref{fig:3d_wing}\textbf{a}, QTT-T is faster: TTI acts on each dimension independently, so only QTTs with $M+1$ cores are rounded. In contrast, QTT-I with TTI requires rounding over $3M$ cores, which is slower. In terms of number of parameters, both formats show exponential compression, although the initial QTT-T compression is better. Moreover, Fig.~\ref{fig:3d_wing}\textbf{a} shows that for a fixed number of sweeps, the time taken by \texttt{TT-Cross} grows approximately linearly with the number of cores and the resulting QTT approximation is inaccurate. Allowing more sweeps improves the \texttt{TT-Cross} quality, but with a polynomial growth in time, making it impractical. In contrast, TTI obtains the QTT representation of the mask in constant time and with a constant error regardless of the number of , see Fig.~\ref{fig:3d_wing}\textbf{a}. Furthermore, Fig.~\ref{fig:3d_wing}\textbf{b} shows \emph{exponential compression} for QTT-I and QTT-T with TTI, since the maximum TT rank remains constant as the number of cores increases. QTT-C achieves slightly better compression because we fix the number of sweeps, and thus the maximum TT-rank remains bounded. However, this compressed tensor exhibits an error that grows exponentially with the number of cores.

\subsection{Synthetic Noise}
\label{Sec: synthetic_noise}

We now turn to synthetic coherent noise, with three-dimensional turbulence as our main example. The mathematical foundations of the constructions used here are reviewed in Appendix~\ref{app: synthetic noise}. In particular, Figs.~\ref{fig:midpoint} and \ref{fig:ptt} show how the \textit{midpoint displacement algorithm} and \textit{Perlin noise} can be reproduced within the TTI formalism. Additional one- and two-dimensional examples are provided in Supplemental Material Sec.~III.

\subsubsection{3D Noise: Turbulence}

As a representative application of synthetic noise, we develop a tensor-network framework for generating three-dimensional turbulence snapshots in our two main architectures, \textit{QTT-interleaved} and \textit{QTT-Tucker}. Although Perlin noise is widely used in procedural modeling \cite{Lagae2010} (see Appendix~\ref{App: perlin} and Fig.~\ref{fig:ptt}), we instead employ a multiscale additive cascade construction (see Appendix~\ref{App:turbulence}). This construction yields a clean Kolmogorov spectrum across the frequency range, something that is difficult to obtain with pure Perlin noise.

We generate a divergence-free velocity field, $\boldsymbol{v}$, encoded as a tensor-network by introducing a vector stream function $\boldsymbol{A}$ such that $\boldsymbol v=\nabla\times\boldsymbol A$. We model $\boldsymbol{A}$'s derivatives via a multiscale cascade (see Eq.~(\ref{eq:turb_ansatz_appendix})). Fixing the lattice size, i.e for a given number of scales $M$, we generate a random QTT-I/QTT-T representing the components of $\boldsymbol{A}$ at each sub-scale $m \in \{1,\ldots,M\}$ and compute their derivatives upsampling them up to the final scale using a cubic B-spline kernel, see see Supplemental Material Sec.~II, together with TTI. The interpolated field's derivatives are rescaled accordingly so that the spectrum follows Kolmogorov's law (see Algorithm~\ref{alg:qtt_Vcascade}).

We collect the statistics over 20 synthetic snapshots of size $2^{10}\times2^{10}\times2^{10}$. Fig.~\ref{fig:synthetic_turbulence}\textbf{a} shows that the snapshots follow the correct Kolmogorov spectrum across the frequency domain (we used two different box sizes to separate QTT-I and QTT-T). Fig.~\ref{fig:synthetic_turbulence}\textbf{b} shows that the synthetic fields exhibit intermittency-like behavior, since the flatness deviates from a 
gaussian profile (kuortosis equal to 3) at small distances.  Moreover, Fig.~\ref{fig:synthetic_turbulence}\textbf{c} shows a linear growth of the maximum rank with respect to the number of cores. Also, we observe that QTT-T scales more favorably, maintaining a low bond dimension and thus resulting in a higher compression.

\section{Discussion and Outlook}
\label{sec: conclusions}

We have presented a general upsampling framework with realizations in both tensor-network (TN) architectures and quantum states. In the TN setting, we introduced \emph{Tensor Train Interpolation} (TTI), a low-rank upsampling scheme that starts from a coarse TN representation and constructs fine-scale cores with controlled TT ranks, while preserving control over the smoothness of the reconstructed signal and guaranteeing a prescribed \(\ell^2\) interpolation error independently of the final resolution. Once this error is fixed, the tail ranks of the resulting QTT remain constant, echoing the fast decay of QTT tail ranks observed in \cite{lindsey_multiscale_2024}; at fixed accuracy, this yields exponential compression. We further showed that our coherent-noise constructions also exhibit low entanglement, extending to this setting the smooth-function behavior analyzed in \cite{haqiu}. Moreover, TTI runs in constant time with respect to the final number of scales. As a result, high-resolution signals in one, two, and three dimensions, ranging from structured functions to procedurally generated fields, can be represented and manipulated with substantially reduced memory and computational cost. Under interleaved encoding, the polynomial TT ranks grow exponentially with the number of spatial dimensions; on the contrary, the tail ranks of QTT-Tucker are dimension independent, leading to faster encoding and higher compression, as observed in \cite{tree_cross}. 
In the quantum setting, we introduced a shallow upsampling circuit that acts on amplitude-encoded functions, or q-samples, and produces smooth approximations whose error scales quadratically with the initial grid spacing, while the circuit depth grows logarithmically with the degree of the quasi-interpolation kernel. The circuit can be implemented using the quantum Fourier transform, a \(\mathrm{poly}(n,m)\) number of controlled one- and two-qubit gates, where \(n\) is the initial number of qubits and \(m\) the number of added qubits, together with an MPS-to-QC encoder.

Beyond compression, data augmentation, and upsampling, our TTI construction opens several promising directions. First, it can provide an accurate low-rank warm start for \texttt{DMRG}-like optimization algorithms. A similar idea was first proposed in~\cite{lubasch_multigrid_2018} and later used in~\cite{coarsetofine} to learn visual data representations of 3D objects. Since these methods operate in spaces whose effective dimension grows exponentially with the number of TT-cores, random initialization becomes increasingly impractical and can lead to poor local minima. By contrast, TTI offers a structured initialization that can substantially improve both robustness and efficiency. This suggests a natural route toward enhanced tensor-network solvers, including PDE pipelines based on variational optimization. Second, TTI produces low-rank approximations that can be mapped directly to shallow quantum circuits through MPS--QC encoders~\cite{Ran2020,Malz2024}. Finally, while our present procedural noise constructions are only weakly band-limited, more refined alternatives such as \emph{Wavelet Noise}~\cite{Cook2005} are nearly perfectly band-limited and appear naturally compatible with our TTI framework.

On the quantum setting, our upsampling quantum circuit can be naturally incorporated into applications based on q-sample states, where an unknown probability distribution is encoded in a quantum state and a higher resolution is required. Under these hypotheses, we can find applications such as fast quantum Monte Carlo estimation~\cite{Montanaro2015}, quantum simulated annealing for combinatorial optimization~\cite{Somma2008}, quantum-walk-based search~\cite{Magniez2011}, speedups for learning agents~\cite{Paparo2014}, Monte Carlo pricing of financial derivatives~\cite{Rebentrost2018}, and financial risk analysis~\cite{WoernerEgger2019}. This last application provides an algorithm that yields a quadratic speedup for risk analysis over Monte Carlo simulations. Finally, a natural direction is to extend the quantum upsampling construction beyond probability-encoded functions to non-negative data encoded directly in quantum amplitudes. This would allow one to use not only quasi-interpolation kernels with quadratic error, but also (quasi-)interpolation kernels with higher convergence rates. Using amplitude encoding rather than probability encoding could also reduce the total circuit depth, since degree-\(p\) polynomial pieces admit exact MPS representations with rank at most \(p+1\), improving the encoding of the polynomial pieces using MPS-to-QC. The main difficulty is that, in this setting, the ancilla register carrying the normalization of the polynomial pieces can no longer be traced out. The upsampling must instead be recovered from conditional measurements of the ancilla and the upsampled state. However, because (quasi-)interpolation kernels decay polynomially, this conditional post-processing requires polynomially many measurements, making it impractical. This motivates the design of measurement-efficient kernels or, alternatively, the development of different interpolation schemes for amplitude-encoded data.

Together, the aforementioned directions suggest that the upsampling of tensor networks and quantum states is not only a practical interpolation or compression tool but also a useful bridge between multiscale scientific computing, synthetic noise, and quantum algorithms.

\section*{Acknowledgments}
We thank Raghavendra Peddinti, Stefano Pisoni, Akshat Shah, Ilia Luchnikov and Giancarlo Camilio for valuable discussions and feedback.

\bibliographystyle{apsrev4-2}  
\bibliography{biblio}  
\clearpage 

\appendix
\makeatletter
\@addtoreset{equation}{section}     
\makeatother
\renewcommand{\theequation}{\Alph{section}\arabic{equation}}

\makeatletter
\@addtoreset{figure}{section}      
\makeatother
\setcounter{figure}{0}
\renewcommand{\thefigure}{\Alph{section}\arabic{figure}}

\section{Tensor Train Representation}
\label{app:TTrep}

In this Appendix, we briefly review the tensor-network formats used throughout this work. We first recall the tensor-train (TT) decomposition, then describe its quantized variant (QTT) for data on dyadic grids, and finally discuss multivariate extensions, including interleaved QTT and QTT-Tucker representations.

\subsection{Tensor Train Representation}

Let \(\mathcal{A}\in\mathbb{R}^{n_{1}\times n_{2}\times \cdots \times n_{d}}\) be a \(d\)-way tensor. In TT form, we introduce three dimensional tensors, called \textit{cores}, $G_{k}\in\mathbb{R}^{r_{k-1}\times n_{k}\times r_{k}},\, k=1,\dots,d,$ with \(r_{0}=r_{d}=1\). Each entry of \(\mathcal{A}\) is written as
\begin{equation}
\label{eq:tt-definition}
\mathcal{A}_{i_{1}\dots i_{d}}
=
G_{1}(i_{1})\,G_{2}(i_{2})\cdots G_{d}(i_{d}),
\end{equation}
where \(G_{k}(i_{k})\in\mathbb{R}^{r_{k-1}\times r_{k}}\) denotes the matrix slice of the \(k\)-th core at index \(i_{k}\). The integers \(r_{1},\dots,r_{d-1}\) are the \emph{TT ranks}, or the dimensions of the bond, and correspond to the ranks of the standard unfolding matrices that separate the indices \((i_{1},\dots,i_{k})\) from \((i_{k+1},\dots,i_{d})\). When these ranks remain small, the TT representation can be exponentially more compact than the full tensor. In the special case \(n_k=2\) for all \(k\), the representation is commonly referred to as a \emph{quantized tensor train} (QTT).

\subsection{Function Encoding in QTT}

For simplicity, let \(f:[0,1)\to\mathbb{R}\) be a function sampled on a uniform dyadic grid of size \(2^{N}\), although the extension to any domain $[a,b)$ is straightforward.  Let the points be uniform distributed $x_{i} \;=\; \frac{i}{2^{N}}, 
\quad
i = 0,1,\dots,2^{N}-1,$
and let \(f_{i} = f(x_{i})\).  Each integer \(i\) admits a binary expansion $i \;=\; a_{1}\,2^{\,N-1} + a_{2}\,2^{\,N-2} + \cdots + a_{N}\,2^{\,0},
\qquad a_{k}\in\{0,1\}.$ Equivalently, let's define $x_{a_{1}\dots a_{N}} 
\;=\; \sum_{k=1}^{N} a_{k}\,2^{-k}$ and $ f_{a_{1}\dots a_{N}} \;=\; f\bigl(x_{a_{1}\dots a_{N}}\bigr).$

Thus the vectorized function values $f(x_i)$ are reshaped into an \(N\)-way tensor of size \(2\times 2\times \cdots \times 2\), whose entries are \(f_{a_{1}\dots a_{N}}\).  A QTT decomposition then writes Eq.~\eqref{eq:tt-definition} with cores of dimensions \(r_{k-1}\times 2 \times r_{k}\) as:
\begin{equation}
\label{eq:qtt-definition}
f_{a_{1}\dots a_{N}}
\;=\; G_{1}(a_{1})\,G_{2}(a_{2})\cdots G_{N}(a_{N}) ,
\end{equation}

When \(f\) is sufficiently smooth, e.g.\ an analytic function or a polynomial of fixed degree, is observed that each the TT-ranks remain small (often independent of \(N\)), so the storage cost $\sum_{k=1}^{N} r_{k-1}\,2\,r_{k} \;=\; \mathcal{O}(N\,r^{2}) $ is exponentially smaller than \(2^{N}\).  In particular, elementary low-rank examples include exponential functions \(e^{\alpha x}\) (rank \(1\)), trigonometric functions \(\sin(\omega x+\phi)\) and \(\cos(\omega x+\phi)\) (rank \(2\)), and any polynomial of degree \(p\) (rank at most \(p+1\)). A rank bound of $\sqrt{\Omega}$ for $\Omega$-bandlimited functions was found in \cite{lindsey_multiscale_2024}, while in \cite{compressibility} several rank bounds were found when the function is replaced by a polynomial expansion.

\subsection{Multivariate Encoding}

There are several ways to encode a $d$-dimensional $2^N\times\cdots\times2^N$ tensor into a QTT-like format: sequential (all bits of each coordinate grouped), interleaved (bit-interleaving across dimensions), or more general tree-tensor networks (see \cite{tindal2024} for different examples), like QTT-Tucker (see Fig.~\ref{fig:TTI}).  While our scheme applies to any, we focus on QTT-interleaved (QTT-I) and QTT-Tucker (QTT-T).

To build the QTT-I format we write each coordinate index $i_m=\sum_{k=1}^N a_{m,k}2^{N-k}$ with bits $a_{m,k}\in\{0,1\}$.  Interleaving weaves bits by significance, yielding a QTT with $dN$ cores 
$G_{m,k}(a_{m,k})\in\mathbb R^{r_{dm+k-m-1}\times r_{dm+k-m}}$ of physical dimension 2.  The functions values are recovered through the contraction
\begin{equation}
\label{eq:interleaved-qtt-cores}
\begin{aligned}
f_{a_{1,1}\dots a_{d,N}}
&=\prod_{k=1}^N\; \left( \prod_{m=1}^d
G_{m,k}\bigl(a_{m,k}\bigr) \right) \\
\end{aligned}
\end{equation}

Alternatively, one may group the $d$ bits at each scale $k$ into a single multi-bit index 
$b_k=(a_{1,k},\dots,a_{d,k})\in\{0,\dots,2^d-1\}$ 
and define cores $H_k(b_k)\in\mathbb R^{r_{k-1}\times r_k}$ of physical dimension $2^d$, so that
\begin{equation}
\label{eq:interleaved-qtt-b}
f_{b_1 \dots b_N}
=H_1(b_1)\,H_2(b_2)\,\cdots\,H_N(b_N).
\end{equation}

Both representations are equivalent, one uses fine grained cores per bit, the other bundled cores per scale, and can be chosen based on implementation convenience.

\subsection{QTT-Tucker encoding}
\label{sec:tucker}

Tucker decomposition \cite{Tucker1966,DeLathauwer2000} provides a multilinear generalization of the singular value decomposition (SVD) to higher-order tensors.  Let 
\(\mathcal{A}\in\mathbb{R}^{n_{1}\times n_{2}\times \cdots \times n_{d}}\) be a \(d\)-way tensor.  
The Tucker model represents \(\mathcal{A}\) as a product of a smaller \emph{core tensor} \(\mathcal{G}\) and a collection of factor matrices $U_{k}\;\in\;\mathbb{R}^{n_{k}\times r_{k}}, \qquad k=1,\dots,d ,$
where each \(r_{k}\) is the dimension of the latent space associated with mode \(k\).  Using the mode-$n$ product \(\times_{n}\),  i.e., for $\mathcal{X}\in\mathbb{R}^{I_1\times\cdots\times I_d}$ and $U\in\mathbb{R}^{J\times I_n}$ one has $(\mathcal{X}\times_n U)_{i_1\ldots i_{n-1}\, j\, i_{n+1}\ldots i_d}=\sum_{i_n=1}^{I_n}\mathcal{X}_{i_1\ldots i_n\ldots i_d}\,U_{j i_n}$, the decomposition is written
\begin{equation}
\label{eq:tucker-definition}
\begin{aligned}
\mathcal{A} &= \mathcal{G}\;\times_{1} U_{1}\;\times_{2} U_{2}\;\cdots\times_{d} U_{d} \\
\mathcal{A}_{i_1 \ldots i_d} &=  \sum_{\gamma_1, \ldots, \gamma_d}\mathcal{G}_{\gamma_1 \ldots \gamma_d} U^{\gamma_1}_{1}\left(i_1\right) \cdots U^{\gamma_d}_{d}\left(i_d\right)
\end{aligned}
\end{equation}
with \(\mathcal{G}\in\mathbb{R}^{r_{1}\times r_{2}\times\cdots\times r_{d}}\).  Each entry of \(\mathcal{A}\) is thus expressed as a multilinear combination of the core entries, modulated by columns of the factor matrices.  The truncated higher-order SVD (HOSVD) computes factor matrices by taking the leading singular vectors of each unfolding, yielding a Tucker representation whose error is quasi-optimal in the Frobenius norm \cite{DeLathauwer2000}.  

The QTT–Tucker format \cite{dolgov2013} or Comb Tensor Networks (CTN) in the physics literature \cite{Chepiga2019}), combines the Tucker and QTT decompositions by representing both the core tensor and the factor matrices in nested low–rank formats (see Fig.~\ref{fig:TTI}).  
Given a Tucker representation $(\mathcal{G},\{U_{k}\})$ of a tensor $\mathcal{A}\in\mathbb{R}^{n_{1}\times\cdots\times n_{d}}$, the core $\mathcal{G}$ is itself decomposed in TT form, $\mathcal{G}_{\gamma_{1}\dots\gamma_{d}}
=G_{1}(\gamma_{1})\,G_{2}(\gamma_{2})\cdots G_{d}(\gamma_{d}),$
with $G_{k}(\gamma_{k})\in\mathbb{R}^{r_{C,k-1}\times r_{C,k}}$, while each factor vector
$U_{k}(:,\gamma_{k})$ is further compressed in QTT representation,
$U(i_{k},\gamma_{k})
=U_{k,1}(\gamma_{k},i_{k,1})\cdots U_{k,L}(i_{k,L}),$ 
after binary encoding the indices $i_{k}\mapsto (i_{k,1},\dots,i_{k,L})$, see Fig.~\ref{fig:ndtt}.  
This two-level structure combines the stability of Tucker with the logarithmic complexity of QTT, leading to efficient storage and computation for high–dimensional data. It further isolates dimensions, allowing efficient operations on each dimension. We note that in 1D QTT and QTT-Tucker are the same, whereas in 2D swapping the order of scales of the first dimension on a sequential QTT gives a QTT-format. Therefore, we would only see remarkable differences in 3 or more dimensions.

\section{General Interpolation Framework}
\label{app: interpolation}

\begin{figure}[t!]
    \centering

    \begin{subfigure}[b]{\linewidth}
        \centering
        \includegraphics[width=1\columnwidth]{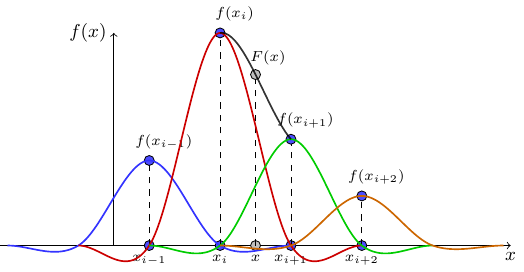}
        \caption{Cubic interpolation}
    \end{subfigure}

    \vspace{1em} 

    \begin{subfigure}[b]{\linewidth}
        \centering
        \includegraphics[width=1\columnwidth]{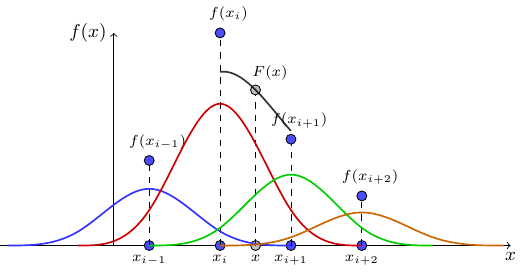}
        \caption{Cubic quasi interpolation}
    \end{subfigure} 
    \caption{\textbf{Interpolation Schemes.} This figure shows different interpolation schemes using 4 data points on a regular grid, $F(x) = \sum_{j=-1}^2
\phi_{3}\!\Bigl(\tfrac{x - x_{j}}{h}\Bigr)\,c_{j} \quad x \in [x_i,x_{i+1})$. Figure (a) shows a $C^1$ cubic interpolant (Keys' kernel \cite{Keys1981}), while figure (b) shows a $C^2$ B-spline quasi interpolant \cite{Schoenberg1988}. The first method has a convergence error rate of $\mathcal{O}(h^3)$, whereas the second has a rate of $\mathcal{O}(h^2)$, where $h$ is the grid size.}
\end{figure}

In this appendix we introduce local polynomial interpolation, its notation, and how it is defined as a kernel convolution for exact interpolation and in a similar way for quasi-interpolation.
Given data $\{(x_{k},\,f(x_{k}))\}_{k=1}^{N}$ on a one‐dimensional grid, we seek a smooth function $F(x)$ that either interpolates or approximates values between gird points.  In full generality we write the interpolation function $F$ as:
\begin{equation}
\label{eq:gen-interp1D}
F(x)
\;=\;
\sum_{\ell=1}^{N} c_{\ell} \, \varphi_{\ell}(x),
\end{equation}
where $\varphi_{\ell}(x)$ are basis (or kernel) functions and $c_{\ell}$ are coefficients chosen by
\begin{itemize}
  \item \emph{Exact interpolation:}  
    $c_{\ell}=f(x_{\ell})$ and 
    $\varphi_{\ell}(x_{j})=\delta_{\ell j}$, so $F(x_{j})=f(x_{j})$.
  \item \emph{Quasi interpolation (e.g.\ B‐splines):}  
    $\{c_{\ell}\}$ can be found by least‐squares or $c_\ell = A_{\ell \ell'}f(x_{\ell'})$ (pre-filtering); in general $F(x_{j})\neq f(x_{j})$.
\end{itemize}

From now on, we will refer to any of the previous cases as just interpolation. In general, each $\varphi_{\ell}(x)$ has compact support in an interval $[x_{\ell-R},\,x_{\ell+R}]$. Therefore, only $2R+1$ terms contribute, allowing for efficient implementation and evaluation.

On a uniform grid $x_k = x_0 + k\,h$, it is common to build
$\varphi_{\ell}(x)$ from a reference piecewise-polynomial function $\phi_{p}$, called \textit{kernel}, of degree $p$,  $\varphi_{\ell}(x) = \phi_{p}\!\Bigl(\tfrac{x - x_{\ell}}{h}\Bigr)$,
giving the local formula
\begin{equation}
\label{eq:local-interp1D}
F(x)
=
\sum_{|\ell - k|\le R}
c_{\ell} \,\phi_{p}\!\Bigl(\tfrac{x - x_{\ell}}{h}\Bigr),
\quad
x\in[x_{k},\,x_{k+1}].
\end{equation}

In general, the interpolation error is $\mathcal{O}(\gamma(f,\phi) h^{p+\delta})$, where $\gamma = \gamma(\phi$, $||f^{(p+1)}(x)||_{\infty})$, and $\delta$ depend on the interpolation scheme and the differentiability class of $f$.

For data on a \(d\)-dimensional Cartesian grid
\(\{\mathbf x_{\mathbf k}=(x_{k_{1}}^{(1)},\dots,x_{k_{d}}^{(d)})\}\)
with values \(f(\mathbf x_{\mathbf k})\), we consider the tensor-product basis
\(\varphi_{\mathbf k}(\mathbf x) = \prod_{m=1}^{d} \varphi_m (x_{m})\),
and define
\begin{equation}
\label{eq:gen-interpD}
F(\mathbf x)
=
\sum_{k_{1}=1}^{N_{1}}\cdots\sum_{k_{d}=1}^{N_{d}}
c_{k_{1}\dots k_{d}} \,\varphi_{\mathbf k}(\mathbf x)\,,
\end{equation}
where, in the interpolation and quasi-interpolation settings, \(c_{k_{1}\dots k_{d}}=f(\mathbf x_{k_{1}\dots k_{d}})\). Compact support in each coordinate guarantees that only \(\prod_{m=1}^{d}(2R_{m}+1)\) neighboring coefficients contribute to the sum at any \(\mathbf x\). Since the kernel \(\varphi_{\mathbf k}(\mathbf x)\) separates across dimensions, the interpolation can be carried out sequentially.

\section{Synthetic Random Fields}
\label{app: synthetic noise}

In this appendix, we review two fundamental techniques for synthesizing random signals and textures: the midpoint displacement algorithm, Sec.~\ref{app: midpoint}, and Perlin noise, Sec.~\ref{App: perlin}. We then describe how these constructions can be extended through fractal superposition to produce continuous, fractal-like fields with scale-dependent roughness.

Midpoint displacement builds fractal profiles through hierarchical subdivision and random perturbations, leading naturally to self-similar structure. In contrast, Perlin noise assigns random gradient vectors to lattice points and interpolates their dot products with local offsets, producing smooth signals with controlled spectral characteristics.

Together with fractal Brownian motion, these methods provide a flexible toolkit for generating synthetic fields with different visual and spectral properties, suitable for terrain modeling, procedural textures, and signal-processing applications.

\subsection{Midpoint Displacement Algorithm}
\label{app: midpoint}
The midpoint displacement algorithm is a recursive method for generating fractal-like terrain profiles.  Starting with two endpoints at positions $x_0$ and $x_N$ with heights $h_0$ and $h_N$, the algorithm repeatedly inserts midpoints, setting each new height to the average of its two neighboring heights plus a random perturbation whose scale decreases with each level of subdivision.

Let $H[0]=h_0$ and $H[N]=h_N$, choose an initial roughness amplitude $R$, and a decay factor $\alpha\in(0,1)$.  The recursion proceeds for levels $\ell=1,2,\dots,k$ where $N=2^k$, as follows:
\begin{enumerate}
  \item At level $\ell$, the segment length is $d=N/2^{\ell-1}$.  For each segment endpoint pair at indices $i$ and $i+d$, compute the midpoint index $m=i+d/2$.
  \item Set $  H[m] = \tfrac12\bigl(H[i]+H[i+d]\bigr) + \mathrm{rand}\bigl([-R,R]\bigr).$

  \item After processing all segments at this level, update $R\gets R\,\alpha$ and proceed to the next level.
\end{enumerate}
The final profile $\{H[0],H[1],\dots,H[N]\}$ exhibits statistical self‑similarity and natural roughness.

\begin{algorithm}[H]
\caption{1D Midpoint Displacement}
\label{alg:midpoint}
\begin{algorithmic}[1]
\Function{MidpointDisplacement}{$N,\;h_0,\;h_N,\;R,\;\alpha$}
  \State Allocate array $H[0\dots N]$, set $H[0]\gets h_0$, $H[N]\gets h_N$
  \For{$\ell=1$ to $k$ where $N=2^k$}
    \State $d \gets N/2^{\ell-1}$
    \For{$i=0, d, 2d, \dots, N-d$}
      \State $m \gets i + d/2$
      \State $H[m] \gets (H[i] + H[i+d])/2 + \mathrm{rand}([-R,R])$
    \EndFor
    \State $R \gets R\,\alpha$
  \EndFor
  \State \Return $H$
\EndFunction
\end{algorithmic}
\end{algorithm}

\begin{figure}[t!]
    \centering
    \includegraphics[width=0.98\columnwidth]{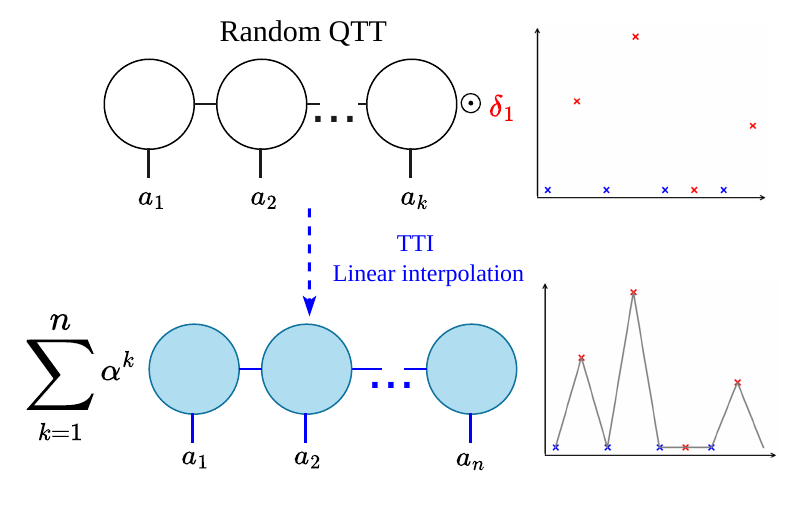}
    \caption{\textbf{QTT midpoint displacement algorithm.} Schematic of a one-dimensional \textit{midpoint displacement algorithm} implemented in \(\mathcal{O}(\log N)\) time and memory using TTI. We begin with a random QTT with \(k\) cores and linearly interpolate it with TTI up to \(n\) scales. We then form a weighted superposition of these linearly extended MPS representations.}
    \label{fig:midpoint}
\end{figure}

In Fig.~\ref{fig:midpoint} we show the transcription of this algorithm in the QTT formalism. Fist we fix the total number of scales $M$, then for each scale $m$ we generate a random QTT, we multiply by a delta function that leaves only odd sites and then this QTT is linearly interpolated $M-m$ more scales. After linear interpolation, we multiply by the roughness factor $\alpha^m$. Finally, the fractal signal is the superposition of the interpolated and rescales random QTTs.

\subsection{Perlin Noise}
\label{App: perlin}

Perlin noise is a procedural texture technique that generates smooth, pseudo-random patterns by interpolating gradient values on a regular lattice.  Classical Perlin uses unit-length gradients, so that each $g_i$ (or $\mathbf g_{\mathbf{i}}$ for multiple dimenisons) has $|g|=1$, yielding $\mathrm{Var}(g)=1$ and a flat gradient spectrum $S_g(\omega)=1$.  Consequently, the noise spectrum is $S_N(\omega)=\lvert H(\omega)\rvert^2,$ where $H(\omega)$ is the Fourier transform of the Perlin fade/interpolation kernel.  In our implementation, we draw gradients with $\mathbb E[g]=0$ and $\mathrm{Var}(g)=1$, so that $S_g(\omega)=1$ and the spectrum remains $S_N(\omega)=\lvert H(\omega)\rvert^2$.  If gradient values become correlated so that $S_g(\omega)\neq1$, then $S_N(\omega)=\lvert H(\omega)\rvert^2\,S_g(\omega),$ allowing controlled “coloring” of the noise by designing the gradient correlation structure.

Perlin noise is often preferred over cubic noise due to its more natural spectral characteristics and isotropy. In the frequency domain, Perlin noise exhibits a smooth, broadband spectrum that decays gradually without introducing artificial cutoffs or ringing artifacts, yielding visually coherent textures across scales. In contrast, cubic noise, constructed from separable tensor‐product kernels, produces anisotropic spectra with rectangular lobes aligned to the coordinate axes, leading to directional artifacts and grid‐aligned features. This inherent isotropy and smoother spectral decay make Perlin noise better suited for generating visually consistent and physically plausible patterns.
Next, we present Perlin noise in one dimension, we exemplify the multidimensional case with three dimensional noise.
\subsubsection{1D Perlin Noise}

At each integer $i\in\mathbb Z$, draw $g_i$ with zero mean and unit variance.  For any $x\in\mathbb R$ set $i_0=\lfloor x\rfloor$, $u=x-i_0\in[0,1)$, and compute $n_0=g_{i_0}\,u$ and $n_1=g_{i_0+1}\,(u-1)$.  Then, choose a fade function $f(t)$, either cubic $f_3(t)=3t^2-2t^3$ for $C^1$ continuity or quintic $f_5(t)=6t^5-15t^4+10t^3$ for $C^2$, and let $s=f(u)$.  Then the 1D Perlin noise is

\begin{equation}
\label{eq:perlin1D}
N(x) = (1 - s)\,n_0 + s\,n_1.
\end{equation}

\begin{algorithm}[H]
\caption{1D Perlin Noise }
\label{alg:perlin1D-compact}
\begin{algorithmic}[1]
\Function{Perlin1D}{$x$}
  \State $i_0\gets\lfloor x\rfloor,\;u\gets x-i_0$
  \For{$a\in\{0,1\}$}
    \State $n_a\gets\textsc{Gradient1D}(i_0+a)\times(u-a)$
  \EndFor
  \State $s\gets f(u)$
  \State \Return $(1-s)\,n_0 + s\,n_1$
\EndFunction
\end{algorithmic}
\end{algorithm}

\subsubsection{3D Perlin Noise}

At each lattice point $(i,j,k)\in\mathbb{Z}^3$, choose a gradient vector
$\mathbf g_{i,j,k}\sim\mathcal{N}(\mathbf0,\mathbf 1)$
(or uniformly on the unit sphere).  For a query point $\mathbf x=(x,y,z)$ set
$i_0=\lfloor x\rfloor$, $j_0=\lfloor y\rfloor$, $k_0=\lfloor z\rfloor$,
and $u=x-i_0$, $v=y-j_0$, $w=z-k_0$.  Apply the fade function
$f_3(t)$ or $f_5(t)$ to get
$s=f(u)$, $t=f(v)$, $r=f(w)$.  Compute the eight corner dot‑products
$n_{abc}
= \mathbf g_{\,i_0+a,j_0+b,k_0+c}\cdot(u-a,\;v-b,\;w-c),
\quad a,b,c\in\{0,1\}.$
Finally, interpolate trilinearly:
\begin{equation}
\label{eq:perlin3D}
\begin{aligned}
n(x,y,z)
&=
(1-r)\Bigl[(1-t)\bigl((1-s)n_{000}+s\,n_{100}\bigr)
 \\& \quad \quad \quad \quad \quad \quad + t\bigl((1-s)n_{010}+s\,n_{110}\bigr)\Bigr]\\
&\quad
+r\Bigl[(1-t)\bigl((1-s)n_{001}+s\,n_{101}\bigr)
   \\
   & \quad\quad\quad\quad\quad\quad+t\bigl((1-s)n_{011}+s\,n_{111}\bigr)\Bigr].
\end{aligned}
\end{equation}
This produces a smoothly varying 3D field whose continuity depends on the choice of $f$. Notice that the previous equation is just the 1D fade function applied to each dimension iteratively.

\begin{algorithm}[H]
\caption{3D Perlin Noise}
\label{alg:perlin3D-compact}
\begin{algorithmic}[1]
\Function{Perlin3D}{$x,y,z$}
  \State $i_0\gets\lfloor x\rfloor,\;j_0\gets\lfloor y\rfloor,\;k_0\gets\lfloor z\rfloor$
  \State $u\gets x-i_0,\;v\gets y-j_0,\;w\gets z-k_0$
  \For{$(a,b,c)\in\{0,1\}^3$}
    \State $n_{abc}\gets\textsc{Gradient3D}(i_0+a,j_0+b,k_0+c)\cdot(u-a,v-b,w-c)$
  \EndFor
  \State $s\gets f(u),\;t\gets f(v),\;r\gets f(w)$
  \For{$c\in\{0,1\}$}
    \State $m_{0c}\gets(1-s)\,n_{0,0,c}+s\,n_{1,0,c},\quad m_{1c}\gets(1-s)\,n_{0,1,c}+s\,n_{1,1,c}$
    \State $p_c\gets(1-t)\,m_{0c}+t\,m_{1c}$
  \EndFor
  \State \Return $(1-r)\,p_0+r\,p_1$
\EndFunction
\end{algorithmic}
\end{algorithm}

\subsection{Fractal Noise}
\label{app:fractal_noise}

\textit{Fractal noise} or \textit{fractal brownian motion} (fbm) takes a smooth noise function $n(\mathbf{x})$ (see Appendices \ref{App: perlin}) as its base and generates a hierarchical superposition of rescaled and attenuated copies, called \emph{octaves}.  The resulting signal exhibits statistical self-similarity across scales and with the proper attenuation it approximates the power-law spectral decay of natural phenomena, such as clouds, fire, water etc.

 Therefore, fractal noise is defined as
\begin{equation}
f(x)
= \sum_{k=0}^{O-1} \alpha^{k}\,
  n\!\bigl(2^{k}\mathbf{x}\bigr),
\label{eq:fractal_noise}
\end{equation}
where $O$ is the number of octaves, and $\alpha\in(0,1)$ controls the amplitude decay (called \emph{persistence}).  The term $n(2^{k}x)$ reproduces finer details at each scale, since it scales the spatial frequency, while $\alpha^{k}$ ensures that the total variance remains finite as $O\!\to\!\infty$.  Therefore, the superposition in~\eqref{eq:fractal_noise} yields a smooth non-periodic signal with fractal characteristics.  Increasing the number of octaves adds finer features, and smaller $\alpha$ results in a faster decay of high-frequency components.  This superposition of signals is fundamental in procedural generation, as it enables controlled roughness together with visually natural complexity.

\subsubsection{QTT Fractal Noise}
\begin{figure}[t!]
    \centering
\includegraphics[width=0.95\columnwidth]{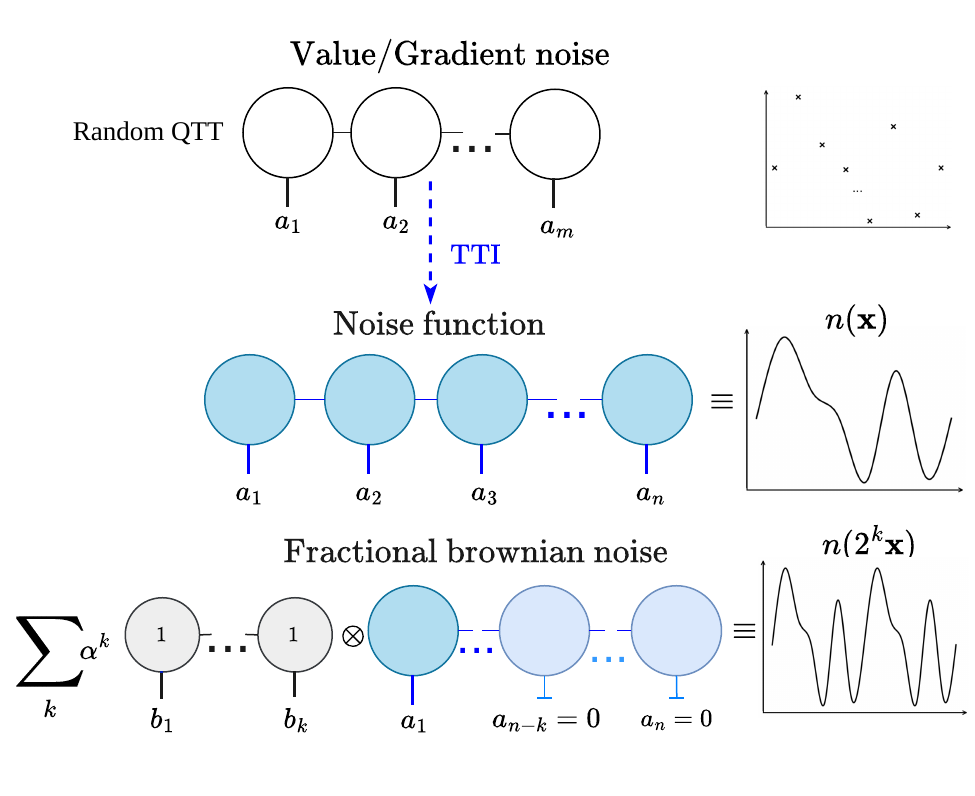}
\caption{\textbf{Fractal noise synthesis.} Schematic of fractal noise construction in QTT format. We begin with a random QTT representing either value noise or gradient noise (Perlin noise), and then apply TTI to obtain a smooth representation. Coherence across scales is straightforward to enforce in QTT format: it is sufficient to evaluate the last cores at \(0\) and append a tensor of ones on the larger scales, with the same physical dimension as the original TT. The final fractal-noise function is obtained as a weighted superposition of rescaled noise functions.}
    \label{fig:ptt}
\end{figure}
Here we describe how to create fractal noise in the TT formalism. Once a noise function $n(\mathbf{x})$ is created as a QTT using TTI, generating octaves is trivial.  For the octave $k$, evaluate the last $k$ cores to $0$ and prepend a QTT of $k$ cores of 1 (bond dimension 1) to create similar copies. This algorithm is illustrated in Fig.~\ref{fig:ptt}.

\section{Synthetic Turbulence}
\label{App:turbulence}

This appendix describes a fully self-contained TT-based algorithm for generating synthetic turbulence with controlled statistics and 
linear bond-dimension scaling.  Our construction enforces incompressibility, reproduces the Kolmogorov spectrum $E(k)\propto k^{-5/3}$, and captures intermittency.

\textit{Multiscale cascade ansatz:} First let's Introduce a vector stream function $\boldsymbol A$ so that $\boldsymbol v=\nabla\times\boldsymbol A$ is divergence-free.  We model its derivatives by a multiscale cascade:
\begin{equation}\label{eq:turb_ansatz_appendix}
\partial_i A_j(x,y,z)
= \sum_{m=2}^{M-1} \omega_m\,\partial_iG^m_{j}(x,y,z),
\end{equation}
where each random TT field $G^m_{j}$ is defined on a $2^m\times2^m\times2^m$ grid, initialized with bond dimension $\chi_G$ and unit variance, and the weights $\omega_m=2^{-4m/3}$ enforce the Kolmogorov scaling.

\textit{TT-based spline interpolation:} To extend each coarse field $G^m_{j}$ to the finest scale without solving global systems, we apply the cubic B-spline quasi-interpolation described in Appendix~\ref{app: interpolation}.  This procedure embeds a TT with $3m$ cores into one with $3M$ cores, maintains $C^2$ smoothness, and preserves linear bond-dimension growth under periodic boundary conditions.

\textit{Derivative evaluation:} Rather than constructing discrete derivative MPOs, based in the 1D derivative operator $D = \frac{1}{2h}(S_1 - S_{-1})$ (which would triple the bond dimension), we differentiate the spline quasi-interpolant directly within the TT format.  In one dimension:
\begin{equation}
    \frac{d}{dx}F(x) =
\sum_{k} f(x_{i+k}) \, \beta'_3\left(\frac{x}{h} -k\right)\, \, \quad x \in [x_i, x_{i+1})
\end{equation}
where $\beta'_3$ is the derivative of the cubic B-spline kernel, see Supplemental Material Sec. IIC.  Expressing $\beta'_3\bigl(x/h-k\bigr)$ as a local quadratic basis yields a TT representation of $F'(x)$ with the same core structure (up to rounding).  This allows computing each $\partial_i A_j$ directly in TT form.

\textit{Velocity reconstruction and bond dimension scaling.} After the interpolation of the derivative, we accumulate the derivatives of the stream function: $\partial_i A_j = \sum_{m=2}^{M-1} \omega_m\,\partial_i G^m_{j},$ and then form the velocity $v_k = \epsilon_{kij}\,\partial_i A_j$ in TT format.  

Algorithm~\ref{alg:qtt_Vcascade} details the full cascade for $\boldsymbol A$ and $\boldsymbol v$.  By construction, the TT bond dimension grows linearly with the number of cores while reproducing incompressibility, the correct energy spectrum, and intermittent fluctuations.

\begin{algorithm}[H]
\caption{3D QTT Velocity Cascade}
\label{alg:qtt_Vcascade}
\begin{algorithmic}[1]
\State $v_x,v_y,v_z \gets \mathbf 0$ \Comment{zero QTT tensors}
\For{$m = 2$ \textbf{to} $M_{\text{scales}}-1$}
    \State $\omega_{m} \gets (2^{-4/3})^{m}$ \Comment{cascade Kolmogorov weight}
    \ForAll{$i,j \in \{x,y,z\}$} \Comment{loop over components}
        \State $G_q^{\text{coarse}} \gets \textsc{RandnQTT}(3m,\chi)$ \Comment{random QTT noise}
        \State $\partial_i G_j^{\text{smooth}} \gets \textsc{D Int}\!\bigl(G_q^{\text{coarse}},M_{\text{scales}}\bigr)$
        \State $v_k \gets \textsc{RoundTT}\!\bigl( v_k + \omega_{m}\epsilon_{ijk} \partial_i \,G_j^{\text{smooth}}\bigr)$
    \EndFor
\EndFor
\State \Return $v_x,v_y,v_z$
\end{algorithmic}
\end{algorithm}

Since the stream vector is expected to belong to the $C^2$ differentiability class, we used a cubic B-spline quasi-interpolant. However, it is also possible to build a $C^2$ interpolant with a quintic kernel without solving a linear system, but that would increase the TT-ranks. Moreover, we also tested cubic interpolants and quadratic quasi-interpolants of class $C^1$. However, we did not observe a significant improvement in the bond dimension. So we kept the cubic implementation instead of the quadratic one used in \cite{Kim2008}.

Even though our noise function is not band limited nor orthogonal in the frequency bands, quadratic/cubic interpolation acts as a low pass filter, making each contribution of Eq.~\eqref{eq:turb_ansatz_appendix} of finite support. Moreover the low-frequency part of the terms with big $m$ are exponentially suppressed by the Kolmogorov scaling $\omega_m$. These two properties of our cascade noise function give rise to the right energy power law $E(k) \propto k^{-5/3}$ (Fig.~\ref{fig:synthetic_turbulence}$\mathbf{a.}$). 

Finally, we can see in Fig.~\ref{fig:synthetic_turbulence}$\mathbf{c.}$ that the synthetic turbulent field has a linear growth with respect to the number of cores. Also, the synthetic turbulent flow exhibits intermittence Fig.~\ref{fig:synthetic_turbulence}$\mathbf{b.}$, something that cannot be achieved with a pure Fourier space construction.

\end{document}


\title{Supplementary Information \\ {\normalsize Efficient upsampling for tensor-network and quantum-state encoded functions}}

\maketitle
%

\section{Supplementary Note 1 - Review of Polynomial Encoding and Shift Matrices}

\vspace{.1cm}

\label{App: plot_shift_encoding}
In this note, we summarize the constructive QTT representation of polynomials given in \cite{oseledets_constructive_2013}. We also explain how to avoid exponentially small entries in the last core. Moreover, we show how to build shift operators as rank-2 MPOs for an arbitrary shift index.

\subsection{Polynomials}
\label{app:qtt_polynomials}

Let \(M(x)=\sum_{k=0}^p c_k x^k\) be a degree-\(p\) polynomial defined on \([0,1)\) and sampled at the points \(x_i=i/2^N\), where \(i=\sum_{m=1}^N a_m 2^{N-m}\), \(a_m\in\{0,1\}\), and \(t_m=a_m 2^{-m}\). Reshaping the vector \(\{M(x_i)\}_{i=0}^{2^N-1}\) into the \(N\)-way tensor \(M_{a_1\cdots a_N}\), we can write its TT decomposition as \(M_{a_1\cdots a_N}=G_1(a_1)\cdots G_N(a_N)\). The cores are

\begin{equation}
\begin{aligned}
G_1(a_1)&=\bigl(\phi_0(t_1),\ldots,\phi_p(t_1)\bigr),\\
[G_k(a_k)]_{i,j}&=
\begin{cases}
\displaystyle\binom{i}{\,i-j\,} t_k^{\,i-j}, & i\ge j,\\
0, & i<j,
\end{cases}\\
G_N(a_N)&=(1,\;t_N,\;\ldots,\;t_N^p)^\top,
\end{aligned}
\end{equation}

where \(\phi_s(x)=\sum_{k=s}^p c_k \binom{k}{s} x^{k-s}\), with \(i,j=0,\dots,p\) and \(2\le k\le N-1\).

We observe that the entries of the last core are proportional to powers of \(2^{-N}\). To avoid this issue, note that this decomposition is based on the binomial expansion \((u+v)^m=\sum_{j=0}^m \binom{m}{j} u^j v^{m-j}\), or, more generally, on polynomial sequences of binomial type, characterized by \(P_n(x+y)=\sum_{k=0}^n \binom{n}{k} P_k(x)P_{n-k}(y)\). Therefore, instead of working directly in the monomial basis, we may first rewrite \(M(x)\) in another polynomial basis and then perform the same construction.

A useful choice is the Abel polynomial sequence, defined by \(A_0(x)=1\) and \(A_n(x)=x(x-an)^{n-1}\) for \(n\ge1\). Choosing \(a=\frac{1-2^N}{2^N}\) avoids exponentially small values in the last core. If we rewrite \(M(x)=\sum_{k=0}^p d_k A_k(x)\), then the corresponding encoding is

\begin{equation}
\begin{aligned}
G_1(a_1)&=\bigl(\phi_0^A(t_1),\ldots,\phi_p^A(t_1)\bigr),\\
[G_k(a_k)]_{i,j}&=
\begin{cases}
\displaystyle\binom{i}{\,i-j\,} A_{i-j}(t_k), & i\ge j,\\
0, & i<j,
\end{cases}\\
G_N(a_N)&=(1,\;A_1(t_N),\;\ldots,\;A_p(t_N))^\top,
\end{aligned}
\end{equation}

where \(\phi_s^A(x)=\sum_{k=s}^p d_k \binom{k}{s} A_{k-s}(x)\), with \(i,j=0,\dots,p\) and \(2\le k\le N-1\).

\subsection{Shift Matrices}
\label{app:shift_matrices}

We now restate the results of~\cite{Kazeev2013}, writing the explicit TT cores of the shift operators. For \(k \in \{0,\dots,2^N-1\}\), we define the shift matrices \(S^{(k)}_{aa'}\), \(R^{(k)}_{aa'}\), and \(L^{(k)}_{aa'}\), with \(a,a' \in \{0,\ldots,2^N-1\}\), by
\begin{equation}
\begin{aligned}
S_{aa'}^{(k)}  &= \delta_{a,(a'-k)\bmod 2^N},\\
R_{aa'}^{(k)}  &= \mathbf{1}_{\{\,a = a' + k < 2^N\,\}},\\
L_{aa'}^{(k)}  &= \mathbf{1}_{\{\,a' = a + k < 2^N\,\}}.
\end{aligned}
\end{equation}
where \(\mathbf{1}_{\{\cdot\}}\) denotes the indicator function, which equals \(1\) if the condition in the subscript is satisfied and \(0\) otherwise. Here, \(S^{(k)}\) is the periodic (circular) shift by \(k\) positions, \(R^{(k)}\) is the right shift by \(k\) positions with zeros at the top, and \(L^{(k)}\) is the left shift by \(k\) positions with zeros at the bottom.

Now write the shift index \(k\) in binary form as \(k=\sum_{i=0}^{N-1}2^i k_i\), and the indices \(a\) and \(a'\) in inverse binary ordering as \(a=\sum_{i=1}^N 2^{N-i} a_i\) and \(a'=\sum_{i=1}^N 2^{N-i} a'_i\), with \(k_i,a_i,a'_i \in \{0,1\}\). Then the shift matrices can be written as rank-2 MPOs:
\begin{equation}
\begin{aligned}
S^{(k)}_{aa'} &= \mathcal{S}_{k_1}(c_1)\, G_{k_2}(c_2)\cdots G_{k_{N-1}}(c_{N-1})\, \mathcal{G}_{k_N}(c_N), \\
R^{(k)}_{aa'} &= \mathcal{R}_{k_1}(c_1)\, G_{k_2}(c_2)\cdots G_{k_{N-1}}(c_{N-1})\, \mathcal{G}_{k_N}(c_N), \\
L^{(k)}_{aa'} &= \mathcal{L}_{k_1}(c_1)\, G_{k_2}(c_2)\cdots G_{k_{N-1}}(c_{N-1})\, \mathcal{G}_{k_N}(c_N),
\end{aligned}
\end{equation}
where \(c_i=(a_i,a'_i)\).

In the equations above, the TT cores take different values depending on the binary digit \(k_i\). Explicitly, they are
\begin{equation}
\begin{aligned}
\mathcal{S}_0(c_i) &=
\begin{bmatrix}
\delta_{a_i a'_i} & \delta_{a_i\,1-a'_i}
\end{bmatrix},
&
\mathcal{S}_1(c_i) &=
\begin{bmatrix}
\delta_{a_i\,1-a'_i} & \delta_{a_i a'_i}
\end{bmatrix}, \\
\mathcal{R}_0(c_i) &=
\begin{bmatrix}
\delta_{a_i a'_i} & \delta_{a_i1}\,\delta_{a'_i0}
\end{bmatrix},
&
\mathcal{R}_1(c_i) &=
\begin{bmatrix}
\delta_{a_i1}\,\delta_{a'_i0} & 0
\end{bmatrix}, \\
\mathcal{L}_0(c_i) &=
\begin{bmatrix}
0 & \delta_{a_i0}\,\delta_{a'_i1}
\end{bmatrix},
&
\mathcal{L}_1(c_i) &=
\begin{bmatrix}
\delta_{a_i0}\,\delta_{a'_i1} & \delta_{a_i a'_i}
\end{bmatrix}, \\
G_0(c_i) &=
\begin{bmatrix}
\delta_{a_i a'_i} & \delta_{a_i1}\,\delta_{a'_i0} \\
0 & \delta_{a_i0}\,\delta_{a'_i1}
\end{bmatrix},
&
G_1(c_i) &=
\begin{bmatrix}
\delta_{a_i1}\,\delta_{a'_i0} & 0 \\
\delta_{a_i0}\,\delta_{a'_i1} & \delta_{a_i a'_i}
\end{bmatrix}, \\
\mathcal{G}_0(c_i) &=
\begin{bmatrix}
\delta_{a_i a'_i} \\ 0
\end{bmatrix},
&
\mathcal{G}_1(c_i) &=
\begin{bmatrix}
\delta_{a_i1}\,\delta_{a'_i0} \\
\delta_{a_i0}\,\delta_{a'_i1}
\end{bmatrix}.
\end{aligned}
\end{equation}

For periodic boundary conditions, \(S^{(k)}\) is used. For more general problems, \(R^{(k)}\) and \(L^{(k)}\) are used together with appropriate boundary conditions.

\section{Supplementary Note 2 - Interpolation Methods}
\label{app:interpolation_methods}

In this note, we collect several techniques for constructing smooth interpolants and approximants from discrete data. Each subsection introduces a different approach, ranging from local polynomial formulas and spline constructions to convolution-based kernels, and provides the basic definitions needed in practice. These methods offer flexibility in balancing smoothness, support width, and computational cost.

\subsection{Local Lagrange Interpolants}
\label{app:lagrange}

A local Lagrange interpolant of degree \(m\) on the interval \([x_i,x_{i+1}]\) is built from a stencil of \(m+1\) nearby nodes containing \(x_i\) and \(x_{i+1}\). Let \(J_i=\{s_i,\dots,s_i+m\}\) denote such a stencil, chosen as centered as possible around the midpoint of \([x_i,x_{i+1}]\). For example, when \(m=3\), a natural centered choice is \(J_i=\{i-1,i,i+1,i+2\}\). Near the boundary, the stencil is shifted so that all indices remain valid.

For \(k\in J_i\), define the local Lagrange basis functions by \(\ell_{i,k}^{(m)}(x)=\prod_{j\in J_i,\;j\neq k}\frac{x-x_j}{x_k-x_j}\), so that \(\ell_{i,k}^{(m)}(x_j)=\delta_{kj}\) for all \(j,k\in J_i\). The interpolant on \([x_i,x_{i+1}]\) is then 

\begin{equation}
    F_i(x)=\sum_{k\in J_i}\ell_{i,k}^{(m)}(x)\,f(x_k).
\end{equation} 
Repeating this construction interval by interval yields a piecewise-polynomial interpolant that uses only local data. Adjacent pieces agree at the grid nodes, so the resulting interpolant is typically \(C^0\), although its derivatives may jump across interval boundaries.

\subsection{Natural Cubic Splines}
\label{app:natural_splines}

A natural cubic spline is the unique piecewise-cubic function \(F\in C^2\) defined on $[x_1,x_N]$ such that \(F(x_k)=f(x_k)\) for all grid points and \(F''(x_1)=F''(x_N)=0\). Thus, unlike the local Lagrange interpolant, the spline is globally \(C^2\).

On a uniform grid, with \(h=x_{k+1}-x_k\) for all \(k\), let \(f_k=f(x_k)\) and \(m_k=F''(x_k)\). The values \(\{m_k\}\) determine the curvature of the spline and satisfy, for \(k=2,\dots,N-1\), the tridiagonal system \(m_{k-1}+4m_k+m_{k+1}=\frac{6}{h^2}(f_{k+1}-2f_k+f_{k-1})\), together with the natural boundary conditions \(m_1=m_N=0\).

Once \(\{m_k\}\) are known, the spline on \([x_k,x_{k+1}]\) is the unique cubic determined by the endpoint values \(f_k,f_{k+1}\) and the endpoint second derivatives \(m_k,m_{k+1}\). Writing \(t=(x-x_k)/h\in[0,1]\), one convenient form is 

\begin{equation}
    F(x)=(1-t)f_k+t\,f_{k+1}
-\frac{h^2}{6}\,t(1-t)\Big((2-t)m_k+(1+t)m_{k+1}\Big).
\end{equation} 
In this form, the first two terms correspond to linear interpolation, while the last two terms add the cubic correction needed to match the curvature and ensure global \(C^2\) smoothness.

Although this linear system can be encoded as a low-rank MPO and solved with a \texttt{DMRG}-like method such as \texttt{ALS} or \texttt{AMEn}, this approach does not by itself provide a theoretical bound on the TT ranks, since the solution MPS is obtained variationally.

\subsection{B--Spline Quasi-Interpolants}
\label{sec:bsplines}

For smooth approximation, one may use the normalized degree-\(p\) B-spline basis \(\{B_k^{(p)}\}\) on a uniform grid of spacing \(h\). There are two equivalent constructions.

\noindent\textbf{Cox--de Boor recursion.}
Define the zeroth-degree basis by \(B_k^{(0)}(x)=\chi(x-x_k)\), where \(\chi(x)=1\) for \(0\le x<h\) and \(\chi(x)=0\) otherwise, i.e. \(\chi\) is the box function of width \(h\). Higher-degree splines are then built recursively as
\begin{equation}
\begin{aligned}
B_k^{(p)}(x)
&=
\frac{x-x_k}{p\,h}\,B_k^{(p-1)}(x)
+
\frac{x_{k+p+1}-x}{p\,h}\,B_{k+1}^{(p-1)}(x),
\qquad p\ge1.
\end{aligned}
\end{equation}

\noindent\textbf{Iterated convolution.}
Equivalently, a degree-\(p\) B-spline is obtained as the \((p+1)\)-fold self-convolution of \(\chi\):
\begin{equation}
B^{(p)}(x)
=
\frac{1}{h^p}
\underbrace{(\chi*\chi*\cdots*\chi)}_{p+1\text{ times}}(x),
\qquad
B_k^{(p)}(x)=B^{(p)}(x-x_k).
\end{equation}
This convolutional view makes the smoothing effect of B-splines explicit.

For quasi-interpolation, it is convenient to use the centered cardinal B-spline kernel \(\beta_p\), obtained by recentering \(B^{(p)}\). On the uniform grid \(x_i=ih\), the simplest convolutional quasi-interpolant is
\begin{equation}
F(x)
=
\sum_{i\in\mathbb{Z}}
f(x_i)\,
\beta_p\!\left(\frac{x-x_i}{h}\right).
\end{equation}
This operator reproduces constants and linear functions. For sufficiently smooth \(f\), it yields a \(C^{p-1}\) approximation with error \(\|f-F\|_{L^\infty}=\mathcal{O}(h^2)\), independently of the spline degree \(p\). Higher-order \(\mathcal{O}(h^{p+1})\) rates are possible with classical spline quasi-interpolants, but then the coefficients are no longer simply \(c_i=f(x_i)\). Instead, they are local linear combinations of nearby samples, i.e. a prefiltering step that produces an effective kernel with larger support.

For the cubic case (\(p=3\)), the centered kernel is
\begin{equation}
\beta_3(t)=
\begin{cases}
\dfrac{1}{6}\bigl(4-6t^2+3|t|^3\bigr), & |t|<1,\\[1mm]
\dfrac{1}{6}\bigl(2-|t|\bigr)^3, & 1\le |t|<2,\\[1mm]
0, & |t|\ge2.
\end{cases}
\end{equation}
Hence,
\begin{equation}
F(x)
=
\sum_{i=1}^{n}
f(x_i)\,
\beta_3\!\left(\frac{x-x_i}{h}\right).
\end{equation}

On each interval \([x_i,x_{i+1})\), setting \(t=(x-x_i)/h\in[0,1)\), the cubic quasi-interpolant can be written in local form as
\begin{equation}
F(x)
=
\sum_{k=-1}^{2}
f(x_{i+k})\,P^{(k)}(t),
\end{equation}
where
\[
P^{(-1)}(t)=\frac{(1-t)^3}{6},\qquad
P^{(0)}(t)=\frac{4-6t^2+3t^3}{6},
\]
\[
P^{(1)}(t)=\frac{1+3t+3t^2-3t^3}{6},\qquad
P^{(2)}(t)=\frac{t^3}{6}.
\]
Equivalently,
\begin{equation}
\left. F(x)\right|_{x\in[x_i,x_{i+1})}
=
\bigl[1\;\;t\;\;t^2\;\;t^3\bigr]\,
M\,
\begin{bmatrix}
f_{i-1}\\
f_i\\
f_{i+1}\\
f_{i+2}
\end{bmatrix},
\end{equation}
with
\[
M=\frac{1}{6}
\begin{bmatrix}
1 & 4 & 1 & 0\\
-3 & 0 & 3 & 0\\
3 & -6 & 3 & 0\\
-1 & 3 & -3 & 1
\end{bmatrix}.
\]

In three dimensions, the tensor-product cubic quasi-interpolant is
\begin{equation}
\begin{aligned}
F(x,y,z)
=
\sum_{i,j,k}
f(x_i,y_j,z_k)\,
\beta_3\!\left(\frac{x-x_i}{h}\right)
\beta_3\!\left(\frac{y-y_j}{h}\right)
\beta_3\!\left(\frac{z-z_k}{h}\right).
\end{aligned}
\end{equation}

Equivalently, one may apply the 1D cubic quasi-interpolant successively along each axis. If \((x,y,z)\in[x_i,x_{i+1})\times[y_j,y_{j+1})\times[z_k,z_{k+1})\), define
\begin{equation}
T_{\alpha,\beta}(z)
=
\mathrm{CA}_z\!\left(
f_{i+\alpha,j+\beta,k-1},
f_{i+\alpha,j+\beta,k},
f_{i+\alpha,j+\beta,k+1},
f_{i+\alpha,j+\beta,k+2};\,z
\right),
\end{equation}
for \(\alpha,\beta\in\{-1,0,1,2\}\), then
\begin{equation}
U_\alpha(y,z)
=
\mathrm{CA}_y\!\left(
T_{\alpha,-1}(z),
T_{\alpha,0}(z),
T_{\alpha,1}(z),
T_{\alpha,2}(z);\,y
\right),
\end{equation}
for \(\alpha\in\{-1,0,1,2\}\), and finally
\begin{equation}
F(x,y,z)
=
\mathrm{CA}_x\!\left(
U_{-1}(y,z),
U_0(y,z),
U_1(y,z),
U_2(y,z);\,x
\right),
\end{equation}
where \(\mathrm{CA}_u\) denotes the 1D cubic approximant along axis \(u\). This stencil-based construction is local, smooth, and avoids the solution of global systems.

\subsection{Cubic Convolution (Keys) Interpolation}
\label{sec:keys}

\begin{figure}[t!]
  \centering
  \includegraphics[width=1\columnwidth]{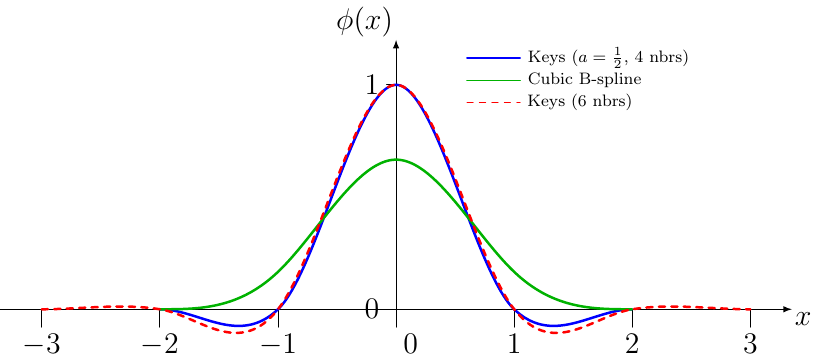}
  \caption{\textbf{Cubic kernels:} a four-neighbor interpolation kernel with \(\mathcal{O}(h^3)\) convergence, a six-neighbor cubic interpolation kernel with \(\mathcal{O}(h^4)\) convergence, and a cubic B-spline quasi-approximation kernel with \(\mathcal{O}(h^2)\).}
  \label{fig:cubic_kernels}
\end{figure}

Cubic convolution, as introduced by Keys~\cite{Keys1981}, provides a \(C^1\) interpolant without solving a linear system for the coefficients, unlike natural splines; see App.~\ref{app:natural_splines}. We start from the kernel-based interpolation formula
\begin{equation}
F(x)=\sum_{k} f(x_k)\,\varphi\!\left(\frac{x-x_k}{h}\right),
\end{equation}
where the Keys kernel is given by (see Fig.~\ref{fig:cubic_kernels})
\begin{equation}
\label{eq:keys}
\varphi(r)=
\begin{cases}
\dfrac{3}{2}|r|^3-\dfrac{5}{2}|r|^2+1, & 0\le |r|<1,\\[0.5ex]
-\dfrac{1}{2}|r|^3+\dfrac{5}{2}|r|^2-4|r|+2, & 1\le |r|<2,\\[0.25ex]
0, & |r|\ge 2.
\end{cases}
\end{equation}
Since \(\varphi\) has support in \([-2,2]\), only four neighboring samples contribute on each interval.

For \(x\in[x_i,x_{i+1})\), set \(t=(x-x_i)/h\in[0,1)\). Then the interpolant depends on the four-point stencil \(\{f_{i-1},f_i,f_{i+1},f_{i+2}\}\) and can be written as
\begin{equation}
F(x)
=
\begin{bmatrix}
1 & t & t^2 & t^3
\end{bmatrix}
M_{\mathrm{Keys}}
\begin{bmatrix}
f_{i-1}\\ f_i\\ f_{i+1}\\ f_{i+2}
\end{bmatrix},
\end{equation}
with
\begin{equation}
M_{\mathrm{Keys}}=
\begin{bmatrix}
0 & 1 & 0 & 0\\[0.25ex]
-\tfrac12 & 0 & \tfrac12 & 0\\[0.25ex]
1 & -\tfrac52 & 2 & -\tfrac12\\[0.25ex]
-\tfrac12 & \tfrac32 & -\tfrac32 & \tfrac12
\end{bmatrix}.
\end{equation}
Equivalently,
\begin{equation}
F(x)
=
\sum_{k=-1}^{2} f(x_{i+k})\,P^{(k)}(t),
\end{equation}
where
\[
P^{(-1)}(t)=-\tfrac12 t+t^2-\tfrac12 t^3,\qquad
P^{(0)}(t)=1-\tfrac52 t^2+\tfrac32 t^3,
\]
\[
P^{(1)}(t)=\tfrac12 t+2t^2-\tfrac32 t^3,\qquad
P^{(2)}(t)=-\tfrac12 t^2+\tfrac12 t^3.
\]

This construction satisfies \(F(x_j)=f(x_j)\) and yields a compactly supported \(C^1\) interpolant.

The extension to higher dimensions is straightforward and is performed sequentially, one dimension at a time. For \((x,y,z)\in[x_i,x_{i+1})\times[y_j,y_{j+1})\times[z_k,z_{k+1})\), define
\begin{equation}
T_{\alpha,\beta}(z)
=
\mathrm{CC}_z\!\left(
f_{i+\alpha,j+\beta,k-1},
f_{i+\alpha,j+\beta,k},
f_{i+\alpha,j+\beta,k+1},
f_{i+\alpha,j+\beta,k+2};\,z
\right),
\end{equation}
for \(\alpha,\beta\in\{-1,0,1,2\}\), then
\begin{equation}
U_\alpha(y,z)
=
\mathrm{CC}_y\!\left(
T_{\alpha,-1}(z),
T_{\alpha,0}(z),
T_{\alpha,1}(z),
T_{\alpha,2}(z);\,y
\right),
\end{equation}
for \(\alpha\in\{-1,0,1,2\}\), and finally
\begin{equation}
F(x,y,z)
=
\mathrm{CC}_x\!\left(
U_{-1}(y,z),
U_0(y,z),
U_1(y,z),
U_2(y,z);\,x
\right),
\end{equation}
where \(\mathrm{CC}_u\) denotes the 1D cubic convolution operator along coordinate \(u\).

The kernel in Eq.~\eqref{eq:keys} has convergence rate \(\mathcal{O}(h^3)\), where \(h\) is the grid spacing. Although it is also possible to construct a \(C^1\) cubic interpolation kernel with \(\mathcal{O}(h^4)\) convergence using a six-point stencil (see Fig.~\ref{fig:cubic_kernels}):
\begin{equation}
\phi(r)=
\begin{cases}
\displaystyle \frac{4|r|^3}{3}-\frac{7|r|^2}{3}+1, & 0\le |r|<1,\\[0.5em]
\displaystyle -\frac{7|r|^3}{12}+3|r|^2-\frac{59|r|}{12}+\frac{5}{2}, & 1\le |r|<2,\\[0.5em]
\displaystyle \frac{|r|^3}{12}-\frac{2|r|^2}{3}+\frac{7|r|}{4}-\frac{3}{2}, & 2\le |r|<3,\\[0.5em]
0, & |r|\ge 3.
\end{cases}
\end{equation}

\section{Supplementary Note 3 - Extra Applications and Metrics}

\vspace{.1cm}

\label{app:extra}
This note  contains additional definitions and metrics that complement the main text. We show one- and two dimensional examples. Since QTT-Tucker in 1D and 2D can be written as a QTT, we only show metrics for QTT-Interleaved, because the matrices for 2D QTT-T are similar.
\subsection{Function Encoding}
This section contains extra examples on function encoding in one- and two- dimensions.

\subsubsection{1D Function}
\label{sec:1dfunciton}
The 1D function used for the metrics is

\begingroup
\small
\setlength{\arraycolsep}{3pt}\renewcommand{\arraystretch}{1.05}
\begin{align}
f(x)
&= A\,\tanh\!\Big(\frac{B(x)+B_g(x)+K_+(x)+K_-(x)}{2.5}\Big), \quad \quad \quad A=1,\nonumber\\[1pt]
B(x)
&= 0.28\,\sin(16\pi x)\,\G{0.20}{0.07} \nonumber\\
&\quad + 0.24\,\cos(44\pi x)\,\G{0.36}{0.05} \nonumber\\
&\quad + 0.20\,\C{5}{18}\,\G{0.58}{0.12} \nonumber\\
&\quad + 0.18\,\sin(120\pi x)\,\G{0.73}{0.03}, \nonumber\\
B_g(x)
&= 0.07\,\sin(2\pi\!\cdot\!1.8\,x+0.2)
 + 0.05\,\cos(2\pi\!\cdot\!3.3\,x+0.9), \nonumber\\
K_+(x)
&= \sum_{(a,s,w)\in\mathcal{H}_+} w\,\pp{x-a}^{3}\,\G{a}{s}, \nonumber\\[-2pt]
\mathcal{H}_+
&= \{(0.22,0.030,{+}0.8),\ (0.37,0.025,{-}0.6),\nonumber\\
&\quad (0.61,0.035,{+}0.7),\ (0.82,0.022,{-}0.5)\}, \nonumber\\
K_-(x)
&= \sum_{(b,s,w)\in\mathcal{H}_-} w\,\pp{b-x}^{3}\,\G{b}{s}, \nonumber\\[-2pt]
\mathcal{H}_-
&= \{(0.28,0.030,{-}0.6),\ (0.42,0.028,{+}0.5),\nonumber\\
&\quad (0.68,0.030,{-}0.5),\ (0.88,0.022,{+}0.4)\}.
\label{eq:fh3}
\end{align}

Notice that this function belongs to $C^2$, since it has discontinuities in the third derivative. 
\endgroup

\begin{figure}[t!]
    \centering
    \includegraphics[width=0.9\columnwidth]{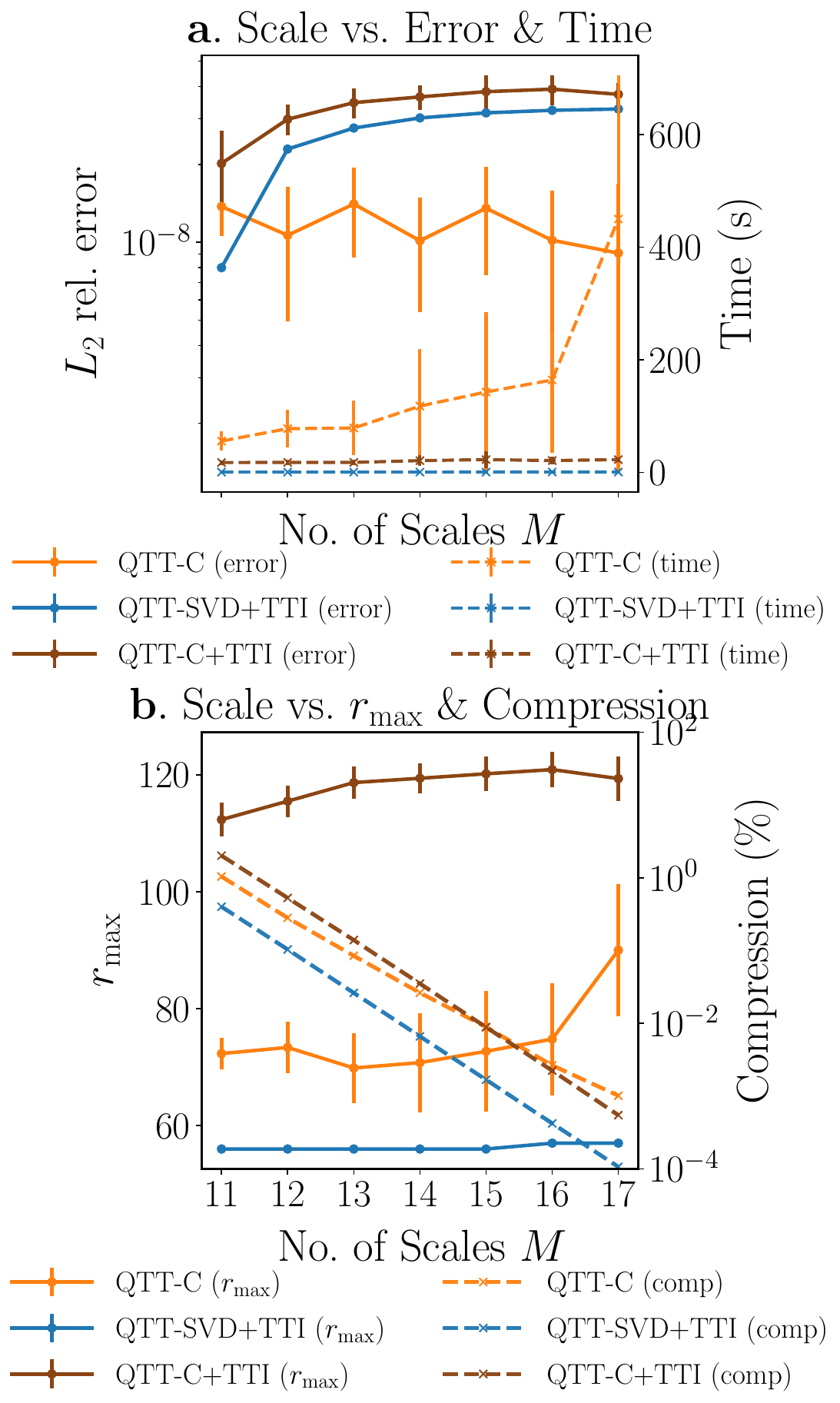}
    \caption{\textbf{QTT 2D correlated gaussian.} This figure compares \texttt{TT-Cross} algorithm vs a coarse \texttt{TT-SVD}/\texttt{TT-Cross} + TTI for a correlated gaussian function in two dimensions using an interleaved encoding. We observe that the time and erank for \texttt{TT-SVD} + TTI outperform \texttt{TT-Cross} plus maintaining a constant error of the same order.}
    \label{fig:corr_gaussian}
\end{figure}

\begin{figure}[htbp!]
    \centering
    \includegraphics[width=0.98\columnwidth]{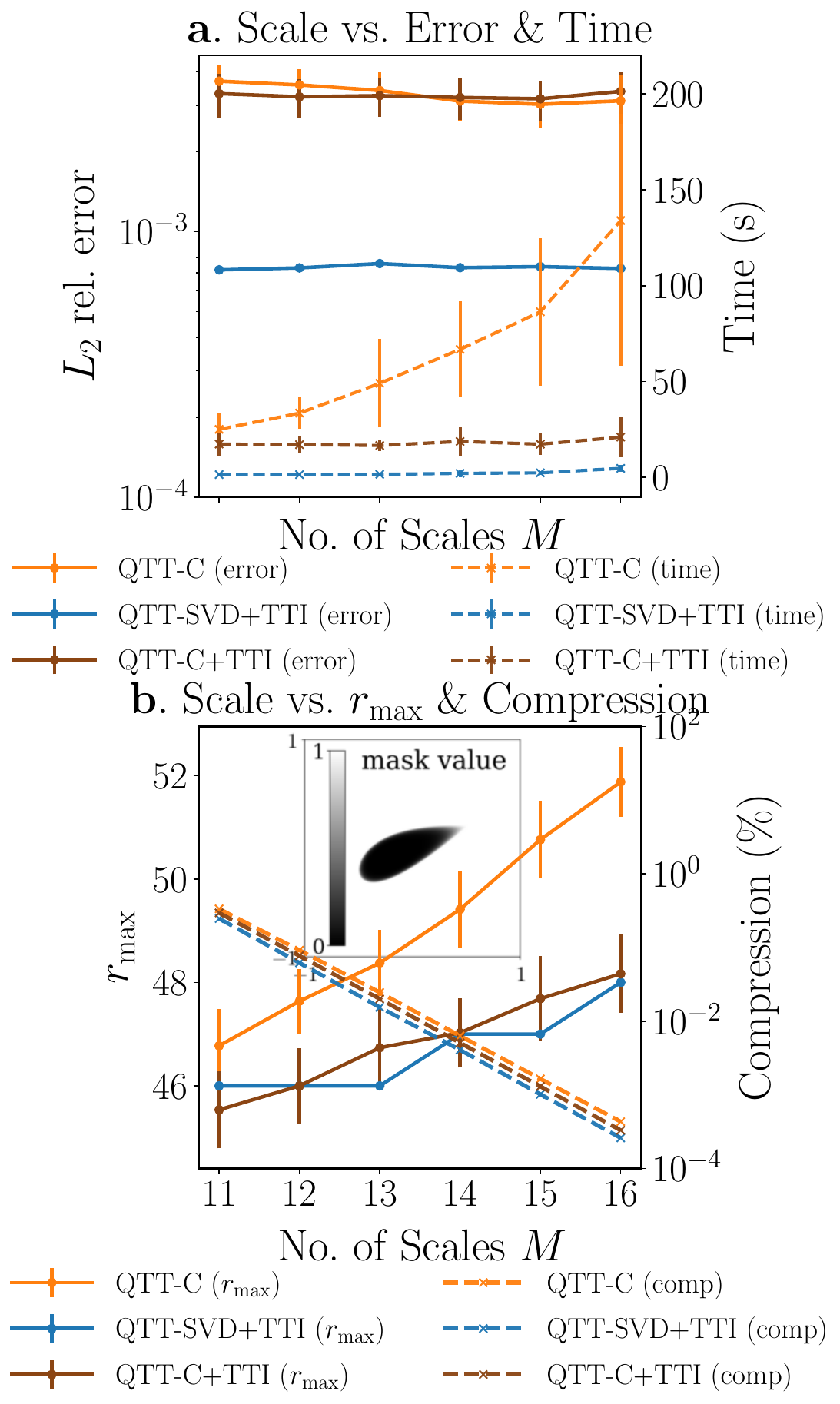}
    \caption{\textbf{QTT 2D Mask.} This figure compares \texttt{TT-Cross} algorithm vs a coarse \texttt{TT-SVD}/\texttt{TT-Cross} + TTI for a softened version of a wing's mask in two dimensions using an interleaved encoding. We observe that the time and erank for \texttt{TT-SVD} + TTI outperform \texttt{TT-Cross}. Moreover, the error found is lower and constant. We observed a linear growing time for full \texttt{TT-Cross} while TTI remains constant. Also, the effective ranks are smaller than the ones found by \texttt{TT-Cross}. }
    \label{fig:wing_mask}
\end{figure}

\subsubsection{Correlated Gaussian Function}

One of many functions that are hard to encode in QTT format is a correlated Gaussian distribution in high dimensions. Since this function is $C^\infty$, we expect an interpolation error of order $\mathcal{O}(h^3)$. In Fig.~\ref{fig:corr_gaussian} we show a constant interpolation error of $10^{-8} \sim h^3$ as expected, both for an initial \texttt{TT-SVD} coarse QTT or a \texttt{TT-Cross} coarse QTT. Moreover, the interpolation time is constant regardless of the final scale. As mentioned in the main text, the first part of the interpolated QTT is bounded by the product of the original TT-ranks times the support of the kernel, then it can be observed in Fig.~\ref{fig:corr_gaussian} that the $\text{rmax}$ for the interpolated QTTs are higher than the \texttt{TT-Cross} one, but since we add low constant tails, the $\text{erank}$ decreases linearly with the number of scales. Therefore, there is a crossover scale, where we get a better compression than \texttt{TT-Cross} in constant time with a fixed error.

\subsubsection{Soft 2D Masks}

As discussed in \cite{NezFernndez2025}, indicator functions do not separate scales; therefore, they are full rank. As proposed in \cite{Peddinti2024}, we can approximate a mask using a $C^\infty$ function. Since the partial derivatives of the softened indicator function are large, we do not expect an error of $\mathcal{O}(h^3)$ at any interpolated scale. Nevertheless, as shown in Fig.~\ref{fig:wing_mask}, our \texttt{TT-SVD} $+$ $\text{TTI}$ finds a better error than \texttt{TT-Cross} in constant time and gives a better compression overall. Therefore, \texttt{TT-SVD} + TTI is a great candidate to encode objects as QTTs to be used for PDEs, TT-based CFD solvers.

\begin{figure}[t!]
    \centering
    \includegraphics[width=0.9\columnwidth]{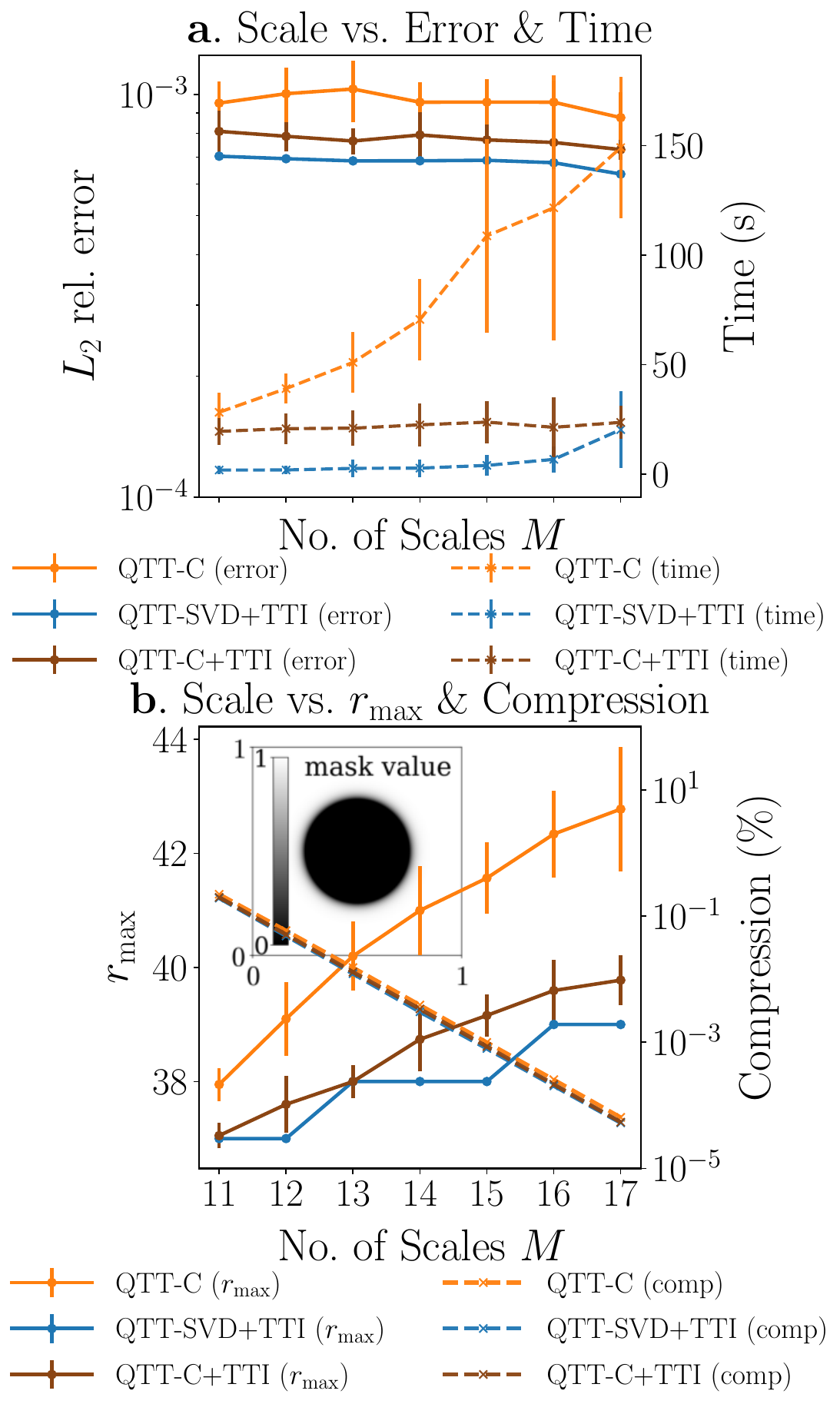}
    \caption{\textbf{QTT 2D Mask.} This figure compares \texttt{TT-Cross} algorithm vs a coarse \texttt{TT-SVD}/\texttt{TT-Cross} + TTI for a softened version of a circle's mask in two dimensions using an interleaved encoding. We observe that the time and erank for \texttt{TT-SVD} + TTI outperform \texttt{TT-Cross}. Moreover, the error found is lower. We observed a linear growing time for full \texttt{TT-Cross} whice TTI remains constant. Also, the effective ranks are smaller than the ones found by \texttt{TT-Cross} for the last grids. }
    \label{fig:circle_mask}
\end{figure}

Moreover, Fig.~\ref{fig:circle_mask} shows the metric for a circle's mask. Since the circle is centered in the middle of the domain, the function can be treated as periodic.

\begin{figure*}[t!]
  \centering
  \setlength{\tabcolsep}{2pt} 

  \begin{subfigure}[b]{0.32\linewidth}
    \includegraphics[width=\linewidth]{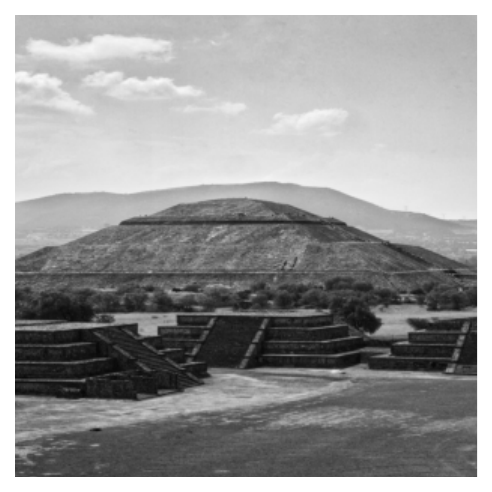}
    \subcaption{Original image ($2048\times 2048$)}
    \label{fig:orig}
  \end{subfigure}\hfill
  \begin{subfigure}[b]{0.32\linewidth}
    \includegraphics[width=\linewidth]{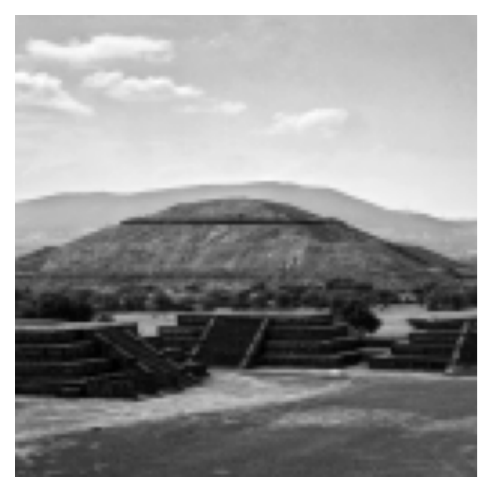}
    \subcaption{Downscaled image ($128\times 128$)}
    \label{fig:down}
  \end{subfigure}\hfill
  \begin{subfigure}[b]{0.32\linewidth}
    \includegraphics[width=\linewidth]{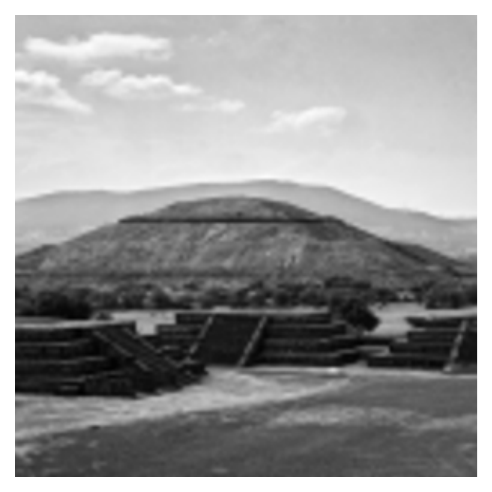}
    \subcaption{Upscaled reconstruction}
    \label{fig:up}
  \end{subfigure}

  \caption{\textbf{Image super-resolution pipeline (visuals).}
  (a) original, (b) downscaled, (c) upscaled reconstruction.
  Quantitative metrics are reported separately in Fig.~\ref{fig:superresolution-metrics}.}
  \label{fig:superresolution}
\end{figure*}

\subsection{Image Super Resolution}
\label{sec: image_superesolution}

\begin{figure}[t!]
  \centering
  \includegraphics[width=0.95\columnwidth]{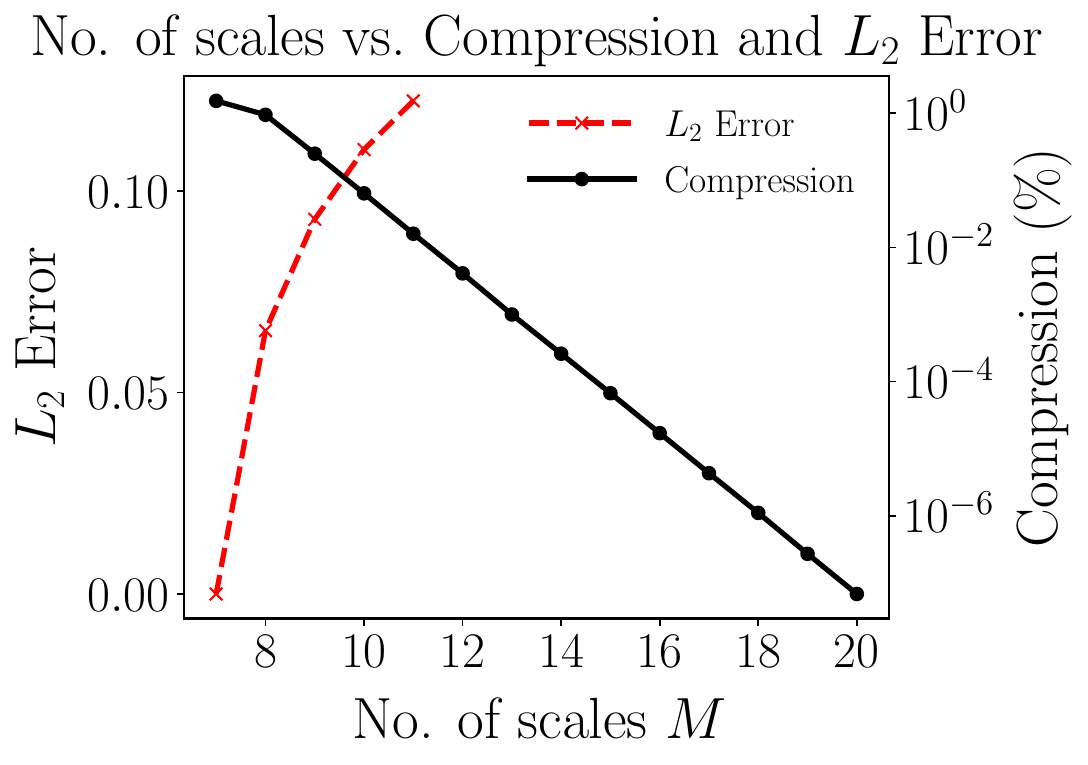}
  \caption{\textbf{Compression and \(\ell^2\) metrics.}
  The compression of the reconstructed image decays \emph{exponentially} with the number of scales.
  We also report the \emph{percentage} \(\ell^2\) error between the upscaled image
  \(\,2^{7}\!\times\!2^{7}\rightarrow 2^{x}\!\times\!2^{x}\,\) and the downscaled image
  \(\,2^{11}\!\times\!2^{11}\rightarrow 2^{x}\!\times\!2^{x}\,\) to match sizes.}
  \label{fig:superresolution-metrics}
\end{figure}

Although AI-based super-resolution methods produce visually compelling results, they often deviate significantly from the original images and can introduce biases inherited from their training data. However, deep learning models, such as CNNs, GANs, transformers, and diffusion networks, now define the state of the art in both perceptual quality and fidelity, bicubic interpolation remains the standard baseline due to its simplicity and reasonable performance. Here, we demonstrate that up-scaled images can be computed and compressed directly in the QTT format, enabling arbitrary-resolution with guaranteed \(\ell^2\) error bounds.

We start with a natural black and white image of
$2^{11} \times 2^{11}$ pixels. Then we down-scaled to $2^7 \times2^7$ pixels, we transform it to QTT-Interleaved using SVD without any compression since at this size SVD is fast and cheap. Then we use TTI in 2D to get the final image. We use Keys cubic interpolation kernel (see Appendix \ref{app:interpolation_methods}), although any Michaell-Netravali filter could be used. Since the last node in our interpolation method is a ghost node, we need to implement QTT boundary conditions manually. Thus we computed the $2^7$ pixel boundaries in QTT format, which again is done with SVD.  
In Fig.~\ref{fig:superresolution} we show the original image, the down-scaled one and the restored image using TTI. We also show the $\ell^2$ error norm together with the compression ratio up to $2^{20} \times 2^{20}$ pixels. For better control of the restored image, e.g. ringing, blur and aliasing, we could have used a custom Mitchell and Netravali filter \cite{Mitchell1988}. Furthermore, since natural images are full rank, the compression shown in Fig.~\ref{fig:superresolution-metrics} is independent of the image, because the initial coarse TT will be full rank.

\subsection{Synthetic Noise}
\label{app: synthetic_noise_extra}

\begin{figure}[t!]
    \centering
    \includegraphics[width=0.92\columnwidth]{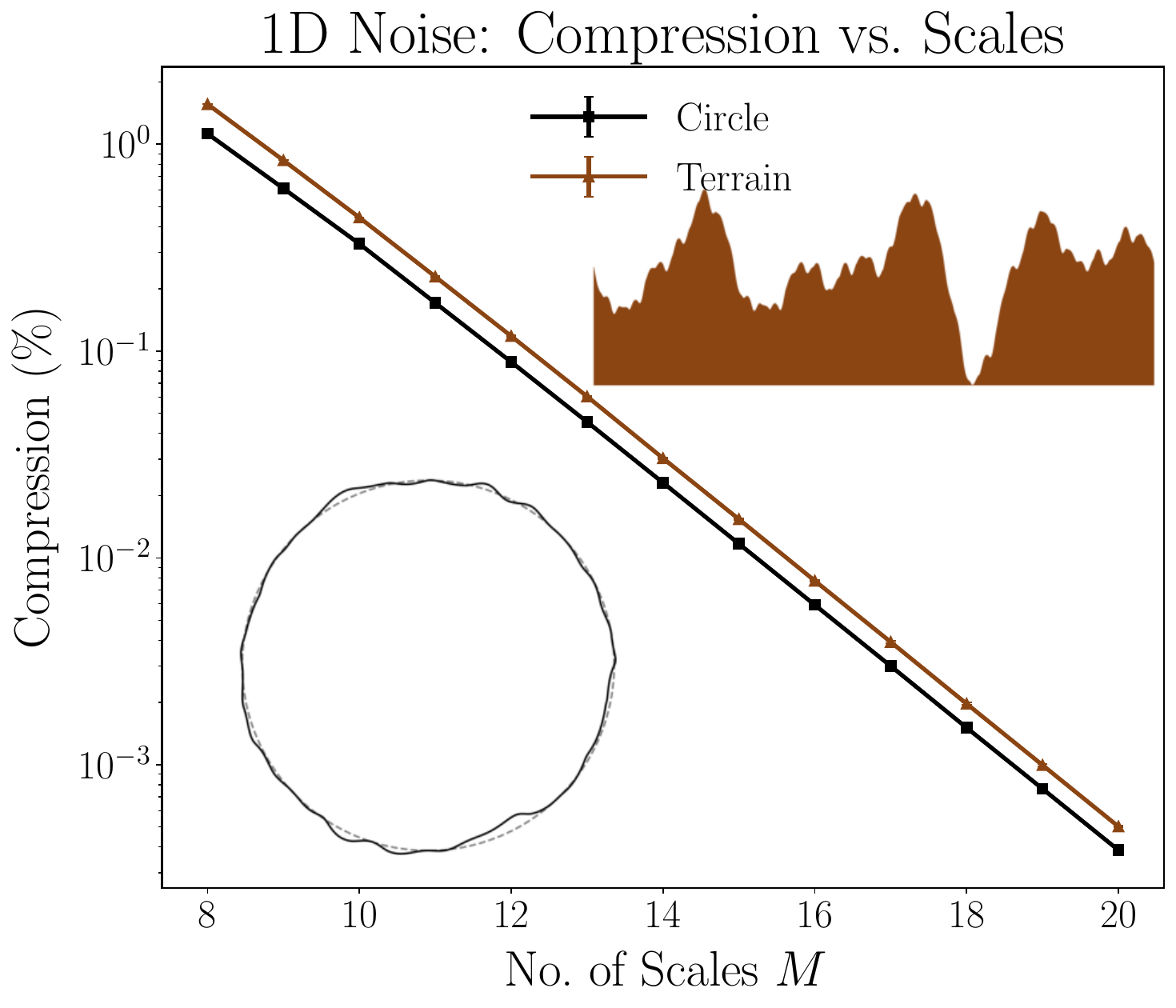}
    \caption{\textbf{1D QTT Perlin noise.} This figure shows two one dimensional applications of fractal Perlin noise. The number of scales $M$ is the total number of cores, then total number of points is $2^M$. On the left we created a \textit{natural} circle, i.e. a circle with embedded noise that makes it look more natural. On the right we show a simple 1D landscape. The circle is made of 3 octaves, whereas the landscape has 4. Both of them have a initial noise grid of $2^3$ gradients. A cubic fade function was used and a persistances of $1/2$ and $1/10$ respectively.}
    \label{fig:1dperlin}
\end{figure}

This section contains examples of one- and two-dimensional applications of synthetic coherent noise. Showcasing Perlin noise and the midpoint algorithm as QTTs.

\subsubsection{1D Noise: Landscapes and offsets}

Here we present a 1D application of coherent noise. In Fig.~\ref{fig:1dperlin} we show a random 1D terrain plotted in 2D. These types of procedural generated terrains are used in side-scrolling games. Another 1D application is presented in Fig.~\ref{fig:1dperlin}, where Perlin noise can be used to create natural hand-written figures.
The total number of points is $N = 2^m$, where $m$ is the number of scales. Both examples can be built with complexity $\mathcal{O}(\log(N))$ and show exponential memory compression.

\subsubsection{2D Noise: Terrains}
\begin{figure}[t!]
    \centering
    \includegraphics[width=0.95\columnwidth]{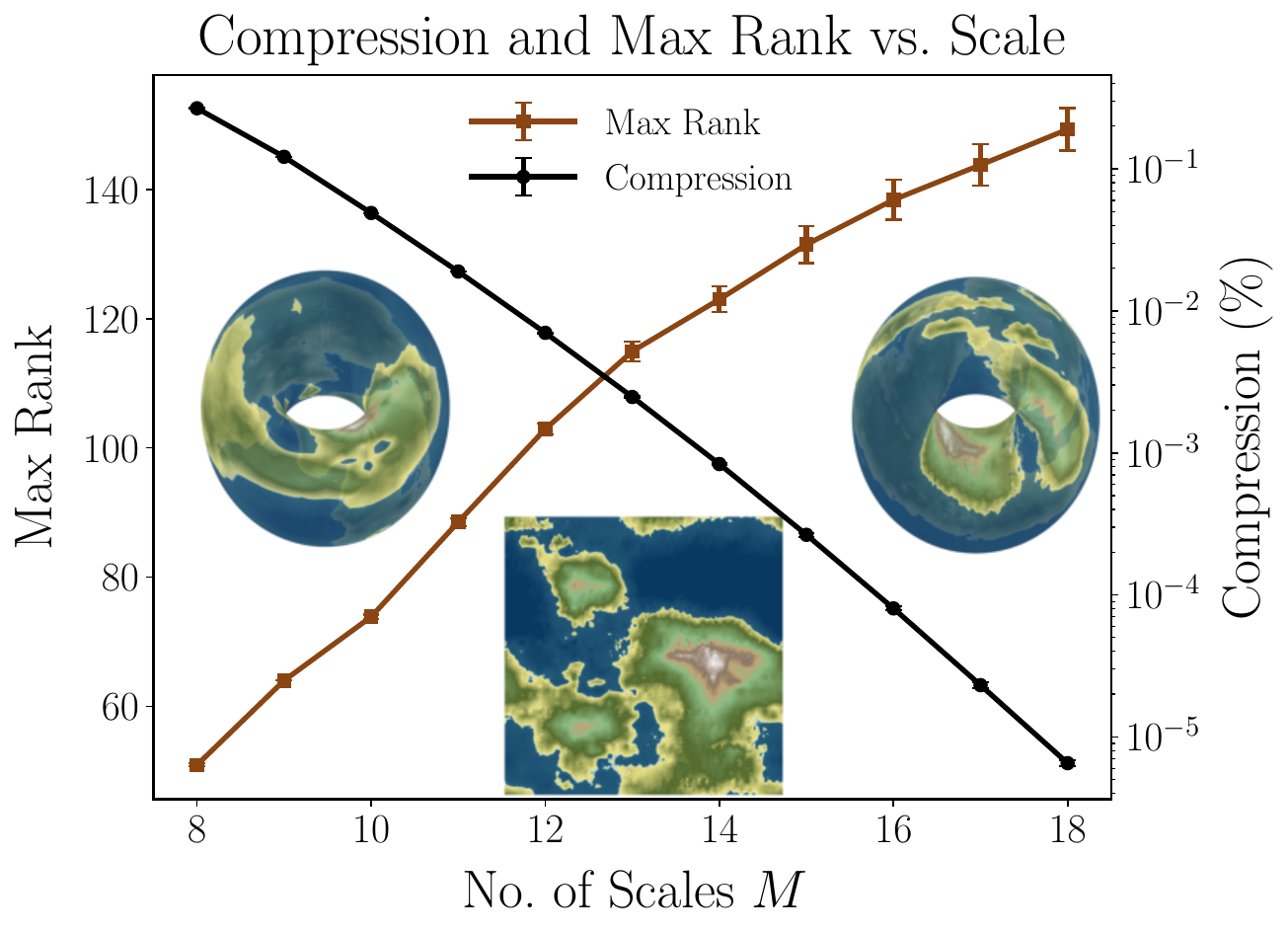}
    \caption{\textbf{2D terrain by midpoint displacement algorithm.} This figure shows the maximum rank and the compression of a 2D terrain generated by our 2D extension of our QTT midpoint displacement algorithm, versus the number of scales $M$, i.e. the number of cores per dimension. We show that the maximum bond dimension grows linearly with the number of scales, achieving an exponential compression as shown by the left line. Each random noise insertion is a QTT of rank 5 and we use the TTI procedure to do bilinear interpolation. We set $\alpha = \frac{1}{2}$ and a rounding of $10^{-6}$, which still captures all the image details. The grid is periodic and the size varies from $2^8 \times 2^8$ up to $2^{18}\times 2^{18}$. We show the colored heatmap of a random generated terrain in 2D and its embedding in a torus.}
    \label{fig:midpoint2d}
\end{figure}

In this two-dimensional application, we demonstrate the generation of random coherent terrains embedded in a three-dimensional space using the midpoint displacement algorithm. Although Perlin or cubic noise methods could be used for this purpose, the midpoint displacement algorithm uses linear interpolation, producing lower-rank QTTs. Although the midpoint algorithm produces some noticeable artifacts, it is the basis for more refined algorithms such as \textit{diamond-square} algorithm \cite{fournier_computer_1982}, which is still used for terrain generation in commercial software. Fig.~\ref{fig:midpoint2d} shows the rank growth together with the compression of randomly generated 2D terrains using the TN version of the midpoint displacement algorithm. 
\vspace{.5cm}

\bibliographystyle{apsrev4-2}
\bibliography{biblio}